\documentclass{article}
\usepackage[a4paper,top=3cm,bottom=2cm,left=3cm,right=3cm,marginparwidth=1.75cm]{geometry}
\usepackage{graphicx}
\usepackage{hyperref}
\usepackage{listings}
\usepackage{amsmath}
\usepackage{amssymb}
\usepackage[toc]{appendix}
\usepackage{pdfpages, caption}
\usepackage{caption}
\numberwithin{equation}{section}
\usepackage[normalem]{ulem}

\newenvironment{texteq}{\refstepcounter{equation}\medskip\par\noindent\begin{minipage}{0.9\textwidth}\begin{footnotesize}}
{\end{footnotesize}\end{minipage}\hfill(\thesection.\arabic{equation})\medskip\par\noindent\ignorespacesafterend}

\newenvironment{texteq*}{\medskip\par\noindent\begin{minipage}{0.9\textwidth}\begin{footnotesize}}
{\end{footnotesize}\end{minipage}\medskip\par\noindent\ignorespacesafterend}

\newcommand{\codett}[1]{\texttt{\small #1}}

\begin{document}

\begin{titlepage}

\begin{center}
{\bf \Large On backpropagating Hessians through ODEs}

\vskip6ex

{\bf Axel Ciceri${}^{(1)}$ and Thomas Fischbacher,${}^{(2)}$   \\ }
\bigskip
${}^{(1)}$
King's College, London\\ 
Department of Mathematics, The Strand, London WC2R 2LS, U.K.
\vskip 5mm
${}^{(2)}$ Google Research\\
Brandschenkestrasse 110, 8002 Z\"urich, Switzerland
\vskip 5mm
\bigskip
\tt{axel.ciceri@kcl.ac.uk, tfish@google.com} \\
\end{center}

\bigskip
\bigskip

\begin{abstract}
We discuss the problem of numerically backpropagating Hessians through
ordinary differential equations (ODEs) in various contexts and
elucidate how different approaches may be favourable in specific
situations.  We discuss both theoretical and pragmatic aspects such
as, respectively, bounds on computational effort and typical impact of
framework overhead.

Focusing on the approach of hand-implemented ODE-backpropagation,
we develop the computation for the Hessian of orbit-nonclosure for a mechanical system.
We also clarify the mathematical framework for extending the backward-ODE-evolution of the
costate-equation to Hessians, in its most generic form.
Some calculations, such as that of the Hessian for orbit non-closure, are performed in a language, defined in terms of a formal grammar, that we introduce to facilitate the tracking of intermediate quantities.

As pedagogical examples, we discuss the Hessian of orbit-nonclosure
for the higher dimensional harmonic oscillator and conceptually
related problems in Newtonian gravitational theory. In particular,
applying our approach to the figure-8 three-body orbit, we readily
rediscover a distorted-figure-8 solution originally described by
Sim\'o~\cite{simo2002dynamical}.

Possible applications may include: improvements to training of `neural ODE'-
type deep learning with second-order methods, numerical analysis of
quantum corrections around classical paths, and, more broadly,
studying options for adjusting an ODE's initial configuration such that
the impact on some given objective function is small.
\end{abstract}

\vfill

\end{titlepage}

\section{Introduction}

The set of dynamical systems that exhibit interesting and relevant
behavior is clearly larger than the set of dynamical systems that can
be explored by relying on symbolic-analytic means only.  Even for
systems where an all-analytic approach is ultimately feasible,
numerics can be a useful tool for exploring their properties and
gaining insights that may guide further exact analysis.

Basic questions that naturally arise in the study of dynamical systems
revolve around finding trajectories that connect~$A$ with~$B$, and
understanding the neighborhood of a given trajectory.  For instance, a
typical task is to demonstrate that a given solution of an ODE indeed
minimizes some particular objective function (rather than merely being
a stationary solution) with respect to arbitrary small variations of
the trajectory. The context in which such questions arise need, of
course, not be related to actual physical motion of some object:
problems such as finding $L_2$-normalizable wave function solutions to
the Schr\"odinger equation in spherical coordinates naturally also
lead to boundary value problems for ODEs.  One popular approach to
solve this is to rely on the `shooting method'\footnote{See
e.g.~\cite{press2007numerical}}.  Clearly, if one can introduce a
smooth measure for how badly off the endpoint of ODE-integration is,
and if one is able to efficiently obtain a numerically good gradient
for this `loss' (to use terminology that is popular in machine
learning), and perhaps even a Hessian, this greatly simplifies finding
viable trajectories via standard numerical optimization.

Overall, the problem of efficiently obtaining good numerical gradients
for (essentially) arbitrary numerical computations that are described,
perhaps not even symbolically, but by any algorithmic implementation,
is solved by sensitivity backgropagation.  The key claim here is that,
if the computation of a scalar-valued function $f$ of a
$D$-dimensional input vector takes effort $E$ and can be performed in
such a way that all intermediate quantities and all control flow
decisions are remembered, then one can transform the code that
evaluates the function into code that efficiently computes the
gradient of $f$, to fill numerical accuracy, at a fixed cost that is a
small multiple of $E$ no larger than eight (see the discussion below),
\emph{irrespective of the number of parameters $D$}!

The basic idea is to first evaluate the (scalar) objective function
normally, but remembering each intermediate quantity used in the
calculation, and allocating a zero-initialized sensitivity-accumulator
for each such intermediate quantity. If the calculation is performed
in such a way that all control flow decisions are remembered, or,
alternatively, all arithmetic operations are recorded on some form of
`tape', one then proceeds backwards through the algorithm, asking for
each intermediate quantity the question of how much the end result
would change if one had interrupted the forward-calculation right
after obtaining the $i$-th intermediate quantity $x_i$ and injected an
ad-hoc $\varepsilon$-change as~$x_i:=x_i+\varepsilon$.  To first order
in $\epsilon$, this would change the result as~$z\to
z+\sigma_i\varepsilon$, where~$\sigma_i$ is the total sensitivity of
the end result with respect to infinitesimal changes to intermediate
quantity $x_i$. Note that, since $\sigma_i$ can be computed from
knowing all internal uses of $x_i$ and the sensitivities of the end
result on all computed-later intermediate quantities that involve
$x_i$, the reason for stepping through the calculation
\emph{backwards} becomes clear. Pragmatically, the procedure is
simplified greatly by considering each individual use of $x_i$ as
providing an increment to the sensitivity-accumulator $\sigma_i$,
since then one can retain the original order of arithmetic statements
rather than having to keep skipping back to the next-earlier
intermediate quantity with unknown sensitivity, then forward to all
its uses. The invention of sensitivity backpropagation is usually
accredited to Seppo Linnainmaa in his 1970 master's
thesis~\cite{linnainmaa1976taylor}, with Bert Speelpenning providing
the first complete fully-algorithmic program transform implementation
in his 1980 PhD thesis~\cite{speelpenning1980compiling}.
Nevertheless, as we will see, some instances of the underlying idea(s)
originated much earlier.

While backpropagation has been rediscovered
also~\cite{rumelhart1986learning} in the narrower context of rather
specific functions that are generally called `neural networks', a
relevant problem that plagued early neural network research was that
demonstrating a tangible benefit from architectures more than three
layers deep was elusive (except for some notable architectures such as
the LSTM~\cite{hochreiter1997long}) until Hinton's 2006 breakthrough
article~\cite{hinton2006fast} that revived the subject as what is now
called the `Deep Learning Revolution'. With the benefit of hindsight,
this work managed to give practitioners more clarity about problems
which, when corrected, mostly render the early `pre-training of deep
networks' techniques redundant.  Subsequent work, in particular on
image classification problems with deep convolutional neural networks,
evolved the `residual neural network' (ResNet) architecture as a
useful design. Here, the main idea is that subsequent layers get to
see the output of earlier layers directly via so-called
`skip-connections', and such approaches have enjoyed a measurable
benefit from high numbers of layers, even in the
hundreds~(e.g.~\cite{he2016deep}).  In the presence of
skip-connections, one can regard the action of each Neural Network
layer as providing (in general, nonlinear) corrections to the output
of an earlier layer, and this view admits a rather intuitive
re-interpretation of a ResNet as a discretized instance of an ordinary
differential equation. With the neural network parameters specifying
the rate-of-change of an ODE, one would call this a `neural
ODE'~\cite{chen2018neural}. This interpretation raises an interesting
point: In case where the ODE is time-reversible\footnote{See
e.g.~\cite{gholami2019anode} for a review on reversing ODEs in a
Machine Learning context.  Our focus here is on time-reversible ODEs,
and the considerations that become relevant when addressing
time-irreversible ODEs, such as they arise in ResNet-inspired
applications, are orthogonal to our constructions.}, we may be able to
do away entirely with the need to remember all intermediate
quantities, since we can (perhaps numerically-approximately)
reconstruct them by running ODE-integration in reverse. Given that
furthermore, the time-horizon over which intermediate quantities get
used is also very limited, one may then also do away with the need for
sensitivity-accumulators.

The procedure of computing sensitivities with backwards-in-time ODE
integration goes as follows. Denote the the forward-ODE as
\begin{equation}
\frac{d}{dt}\vec y=\vec F(\vec y)\,,
\end{equation}
(where the time-dependency on the state-vector $\vec y(t)$ has been suppressed in the notation),
with initial state at $t=t_S$ and end state at $t=t_E$ given, respectively, by the~$D$-vectors~$\vec y_S$ and~$\vec y_E$.
Consider then an endpoint-dependent scalar `loss' $L(\vec y_E)$.
The recipe to obtain the sensitivity of~$L$ with respect to the coordinates of
$\vec y_S$ (i.e.~the gradient of the loss as a function of $\vec y_S$)
is to ODE-propagate from $t=t_E$ back to $t=t_S$ the doubled-in-size
state vector
\begin{equation}
\begin{split}
\vec Y_{t=t_E}=(&y_{E,0},\,y_{E,1},\, \ldots, \,y_{E,D-1};\\
 &\partial L/\partial y_{E, 0},\, \partial L/\partial y_{E, 1}, \,\ldots, \partial L/\partial y_{E, D-1})
\end{split}
\end{equation}
according to the doubled-in-state-vector-size ODE\footnote{
We use the widespread convention that a comma in an index-list
indicates that the subsequent indices correspond to partial
derivatives with respect to the single positional dependency
of the quantity, so: $(\partial /\partial x_i f_k(\vec x)=:f_{i,k}(\vec x)$.}
\begin{equation}
  \label{eq:ode2}
\begin{split}
 \frac{d}{dt}Y_{k, 0\le k<D} &= F_k(Y_0, Y_1, \ldots, Y_{D-1})\\
 \frac{d}{dt} Y_{D+k, 0\le k<D} &= -Y_{D+i} F_{i, k}(Y_0, Y_1, \ldots, Y_{D-1})\,.
\end{split}
\end{equation}

Here, the first half of the state-vector reproduces intermediate
states, and the $k$-th entry of the second half of the state vector
answers, at time $t$, the question by how much the final loss $L$
would change relative to $\varepsilon$ if one had interrupted the
forward-ODE-integration at that time $t$ and injected an $\varepsilon$
jump for state-vector coordinate $k$.\footnote{One is readily
convinced of this of this extended ODE by considering regular
sensitivity backpropagation through a simplistic time integrator that
uses repeated $\vec y_{n+1} = \vec y_{n} + \delta T\cdot F(\vec y_n)$
updates with small $\Delta T$. The form of the result cannot depend on
whether a simplistic or sophisticated ODE integrator was used as a
reasoning tool to obtain it.}

This expanded ODE, which addresses the general problem of
ODE-backpropagating sensitivities, has historically arisen in other
contexts.  Lev Pontyragin found it in 1956 when trying to maximize the
terminal velocity of a rocket (see~\cite{gamkrelidze1999discovery} for
a review) and so arguably may be regarded as an inventor of
backpropagation even earlier than Linnainmaa. However, one finds that
the second (``costate'') part of this expanded ODE is also closely
related to a reinterpretation of ``momentum'' in Hamilton-Jacobi
mechanics (via Maupertuis's principle), so one might arguably consider
Hamilton an even earlier originator of the idea of
backpropagation\footnote{The underlying connections will be explained
in more detail elsewhere~\cite{RosettaStone}. One might want to make
the point that, given the ubiquity of these ideas across ML, Physics,
Operations Research, and other disciplines, this topic should be made
part of the standard curriculum with similar relevance as linear
algebra.}.

Pragmatically, a relevant aspect of Eq.~(\ref{eq:ode2}) is that, since
we are reconstructing intermediate quantities via ODE-integration, we
do also not have the intermediate quantities available that did arise
during the evaluation of $F$.
In some situations, the Jacobian~$F_{i,k}\equiv(\partial/\partial y_k)F_i$,
which has to be evaluated at position~$\vec y(t)$, has a simple analytic structure that one would
want to use for a direct implementation in hand-backpropagated code.
If that is not the case, one should remind oneself that what is asked
for here is not to first compute the full Jacobian, and then do a
matrix-vector multiplication with the costate:
doing so would involve wasting unnecessary effort, since the computation of the full Jacobian
turns out to not being needed. Rather, if $\vec\sigma(t)$ is the
sensitivity of the end result on $\vec y(t)$, and $\vec y(t+\Delta
t)=\vec y(t)+\Delta t\cdot \vec F(\vec y(t))$, then computing
$F_{i,k}\sigma_i$ amounts to computing the contribution of the
sensitivity of the end result on $\vec y(t)$ that is coming from the
$\Delta t\cdot \vec F(\vec y(t))$ term. So, one is asked to
sensitivity-backpropagate the gradient into the argument of
$F(\vec y(t))$. Since none of the intermediate quantities from the original
forward-ODE-integration calculation of this subexpression have been
remembered, we hence have to redo the computation of $F(\vec y(t))$ on
the backward-pass, but in such a way that we remember intermediate
quantities and can backpropagate this subexpression.
Given that the difference between doing a matrix-vector product
that uses the full Jacobian and backpropagating an explicit
sensitivity-vector through the velocity-function might be unfamiliar
to some readers, we discuss this in detail in
Appendix~\ref{app:velocity_jacobian}.

The impact of this is that we have to run a costly step of forward
ODE-integration, namely the evaluation of $F(\vec y)$, one extra time
in comparison to a backpropagation implementation that did remember
intermediate results. So, if the guaranteed maximal ratio
\emph{\{gradient calculation effort\}}:\emph{\{forward calculation
effort\}} for working out a gradient by regular backpropagation of a
fixed-time-step ODE integrator such as RK45 with all
(intermediate-time) intermediate quantities were remembered were $b$,
having to do an extra evaluation of $F(\vec y)$ on the backward step
makes this $B:=b+1$ if intermediate quantities have to be
reconstructed.

Knowing how to backpropagate gradients through ODEs in a
computationally efficient way (in particular, avoiding the need to
remember intermediate quantities!) for ODEs that admit reversing
time-integration certainly is useful. Obviously, one would then in
some situations also want to consider higher derivatives,
in particular Hessians. One reason for this might be to use (perhaps
estimated) Hessians for speeding up numerical optimization
with 2nd order techniques, roughly along the lines of the BFGS
method.
In the context of neural networks,
this has been explored in~\cite{liu2021second}.
Other obvious reasons might include trying to understand stability
of an equilibrium (via absence of negative eigenvalues in the Hessian),
or finding \hbox{(near-)degeneracies} or perhaps \hbox{(approximate-)extra} symmetries
of trajectories, indicated by numerically small eigenvalues of the Hessian.

In the rest of this article, when a need for an intuitive pedagogical
and somewhat nontrivial example arises, we will in general study
orbit-nonclosure of dynamical systems where the state of motion is
described by knowing the positions and velocities of one or multiple
point-masses. The `loss' then is coordinate-space squared-distance
between initial-time and final-time motion-state for a given time-interval.

\section{Options for propagating Hessians through ODEs}

In this Section, we comprehensively discuss and characterize the major
options for propagating Hessians through ODEs.

\subsection{Option 1: Finite differencing}
\label{sec:optionFD}
An obvious approach to numerically estimating a Hessian is to run
ODE-integration at least $D^2+D+1$ times to determine the
corresponding number of parameters in the function's value, gradient,
and Hessian via finite differencing. One of the important drawbacks of
this approach is that finite differencing is inherently limited in its
ability to give good quality numerical derivatives. Nevertheless, the
two main aspects that make this approach relevant is that (a) it does
not require the ODE to admit numerical reverse-time-integration, and
(b) it is easy to implement and so provides an easy-to-understand
tractable way to implement test cases for validation of more
sophisticated numerical methods (which in the authors' view should be
considered mandatory). Figure~\ref{fig:fdpy} shows the generic
approach to finite-differencing, using Python\footnote{While some
effort related to unnecessary allocation of intermediate quantities
could be eliminated here, this is besides the point for an
implementation intended for providing (approximate) ground truth for
test code. The default step sizes have been chosen to give good
numerical approximations for functions where all values and
derivatives are numerically `at the scale of 1', splitting available
IEEE-754 binary64 floating point numerical accuracy evenly for
gradients, and in such a way that the impact of 3rd order corrections
is at the numerical noise threshold for Hessians.}.

\begin{figure}
{\footnotesize
\begin{lstlisting}[language=Python]
def fd_grad(f, x0, *, eps=1e-7):
  """Computes a gradient via finite-differencing."""
  x0 = numpy.asarray(x0)
  if x0.ndim != 1:
    raise ValueError(f'Need 1-index position-vector x0, got shape: {x0.shape}')
  x0_eps_type = type(x0[0] + eps)
  if not isinstance(x0[0], x0_eps_type):
    # If `eps` changes cannot be represented alongside x0-coordinates,
    # adjust the array to have suitable element-type.
    x0 = x0.astype(x0_eps_type)
  dim = x0.size
  f0 = numpy.asarray(f(x0))
  result = numpy.zeros(f0.shape + (dim,), dtype=f0.dtype)
  denominator = 2 * eps
  xpos = numpy.array(x0)
  for num_coord in range(dim):
    xpos[num_coord] = x0[num_coord] + eps
    f_plus = numpy.asarray(f(xpos))
    xpos[num_coord] = x0[num_coord] - eps
    f_minus = numpy.asarray(f(xpos))
    result[..., num_coord] = (f_plus - f_minus) / denominator
    xpos[num_coord] = x0[num_coord]
  return result

def fd_hessian(f, x0, *, eps=1e-5):
  grad_f = lambda x: fd_grad(f, x, eps=eps)
  return fd_grad(grad_f, x0, eps=eps)

def verified_grad(f, fprime, eps=1e-7, rtol=0.1, atol=1e-3):
  def get_and_compare_grad(x0):
    grad_via_fprime = fprime(x0)
    grad_via_fd = fd_grad(f, x0, eps=eps)
    if not numpy.allclose(grad_via_fprime, grad_via_fd, rtol=rtol, atol=atol):
      raise ValueError('Gradient mismatch!')
    return grad_via_fprime
  return get_and_compare_grad  
\end{lstlisting}
}
\begin{caption}
  {Basic finite differencing for a multi-index array, using NumPy.
    The `\codett{verified\_grad}' wrapper is useful during code development
    to automatically assert that an algorithmic gradient is aligned with
    the finite-differencing gradient across all individual components.
    A computationally less expensive smoke test might want to check this
    only along one random direction only.}
\label{fig:fdpy}
\end{caption}
\end{figure}

\subsection{Option 2: Autogenerated (backpropagation)${}^2$}

In general, the design philosophy underlying Machine Learning
frameworks such as TensorFlow~\cite{tensorflow2015-whitepaper},
PyTorch~\cite{paszke2019pytorch}, JAX~\cite{jax2018github},
ApacheMxnet~\cite{chen2015mxnet}, and
also~DiffTaichi\cite{hu2019difftaichi} emphasizes the need to remember
intermediate quantities. With this, they can synthesize gradient code
given a computation's specification automatically.  These frameworks
in general have not yet evolved to a point where one could use a
well-established approach to exercise a level of control over the
tensor-arithmetic graph that would admit replacing the need to
remember an intermediate quantity with a (numerically valid)
simulacrum in a prescribed way.  As such, grafting calculations that
do not remember intermediate quantities is generally feasible (such as
via implementing low level TensorFlow op in C++), but not at all
straightforward. The article~\cite{chen2018neural} shows
autograph-based code that can in principle stack backpropagation to
obtain higher derivatives and while this option is readily available,
there are two general problems that arise here:

\begin{itemize}

\item This option is only available in situations where one can deploy
  such a computational framework. This may exclude
  embedded applications (unless self-contained low level code that
  does not need a ML library at run time can be synthesized from a
  tensor arithmetic graph), or computations done in some other
  framework (such as some popular symbolic manipulation packages like
  Mathematica~\cite{wolfram2003mathematica} or Maple~\cite{heck1993introduction}) that cannot
  reasonably be expected to see generic backpropagation capabilities
  retrofitted into the underlying general-purpose programming language
  (as for almost every other complex programming language).

\item The second backpropagation will use the endpoint of the first
  ODE-backpropagation as the first half of its extended state-vector.
  This means that the first quarter of this vector will be the initial
  position, but reconstructed by first going through $t_S\to t_E$
  ODE-integration and then backwards through $t_E\to t_S$, accruing
  numerical noise on the way. In general, one would likely prefer to
  use the noise-free initial position as part of the starting point of
  the ODE-backpropagation of the first ODE-backpropagation, to keep
  the numerical error of the end result small.\footnote{Some
  experimental data about this is shown in~\cite{liu2021second}.}
\end{itemize}

The obvious advantage of this approach is of course convenience --
where it is available, one need not concern oneself with writing the
code and possibly introducing major mathematical mistakes.

\subsection{Option 3: Hand-implemented (backpropagation)${}^{2}$}
\label{sec:handbackprop}

The most direct way to address the problem that framework-generated
backpropagation code will in general use an unnecessarily noisy
reconstructed initial position, and also to backpropagate Hessians in
situations where no such framework is available, is to write the code
by hand. In practice, the exercise of implementing custom code for a
backpropagation problem can easily become convoluted.  This is due, in
part, to the accumulation of inter-dependent quantities which
accumulate through the various intermediate steps of the computation.
In an effort to address this issue, we show in
Section~\ref{sec:handcrafted}, by means of a nontrivial example, how
the emerging complexity can be managed by using some dedicated
formalism to precisely describe dependencies of tensor-arithmetic
expressions. The formalism itself and its underlying design
principles are described in Section~\ref{app:tcd_dsl}.  In the present
Subsection, we restrict to an overview of the general methodologies
and discuss their costs.

Since backpropagation in general considers sensitivity of a scalar
result on some intermediate quantity, we have a choice here to either
consider backpropagation of every single entry of the
sensitivity-gradient to build the Hessian row-by-row, using a
state-vector of length $4D$ for the
backpropagation-of-the-backpropagation, where $D$ is the state-space
dimension of the problem, or considering independent sensitivities
with respect to each individual entry of the gradient.

Since the basic reasoning behind sensitivity backpropagation also
applies here, at least when considering a non-adaptive (i.e. fixed
time step) ODE integration scheme, such as RK45, the total effort to
work out a Hessian row-by-row is upper-bounded by $B^2DE$, where $E$
is the effort for the forward-computation and $B$ is the small and
problem-size-independent multiplicative factor for using sensitivity
backpropagation we discussed earlier in this work. So, the total
effort for obtaining the gradient is $\le BE$. Hence, the effort to
obtain one row of the Hessian is $\le B^2E$.

Computing the Hessian by performing a single
backpropagation-of-the-backpropagation ODE integration that keeps
track of all sensitivities of individual gradient-components on
intermediate quantities in one go correspondingly requires the same
effort\footnote{This claim holds as long as ODE integration does not
try to utilize Jacobians. This is the case (for example) for
Runge-Kutta type integration schemes, but not for numerical ODE
integration schemes such as CVODE's CVDense,
cf.~\cite{cohen1996cvode}.}, but
trades intermediate memory effort for keeping track of state-vectors
that is $\propto D$ for $2D$ independent ODE-integrations
(backpropagating the ODE-backpropagation, and then also the original
ODE-integration) for intermediate memory effort that is $\propto D^2$
for two ODE-integrations.

In (arguably, rare) situations where the coordinate basis is not
randomly-oriented but individual components of the ODE have very
specific meaning, it might happen that the ODEs that involve $F_{m,i}$
for different choices of index $i$ become numerically challenging to
integrate at different times. In that case, when using an adaptive ODE
solver, one would expect that the time intervals over which
integration has to proceed in small steps is the union of the time
intervals where \emph{any} of the potentially difficult components
encounters a problem. In such a situation, it naturally would be
beneficial to perform the computation of the Hessian row-by-row. 
As such, it might also make sense to look for a coordinate basis
in which different parts of the Hessian-computation get difficult
at different times, and then perform the calculation row-by-row.

\subsection{Option 4: FM-AD plus RM-AD}

In general, there is the option to obtain Hessians is by combining
forward mode automatic differentiation (FM-AD) with reverse mode
automatic differentiation (RM-AD), i.e. backpropagation, where each
mode handles one of the two derivatives.  In general, FM-AD gives
efficient gradients for $\mathbb{R}^1\to \mathbb{R}^D$ functions, and
RM-AD gives efficient gradients for $\mathbb{R}^D\to \mathbb{R}^1$
functions, with effort scaling proportionally if one goes from the
one-dimensional input / output space to more dimensions.

If FM-AD is used to obtain the gradient, and then RM-AD to get the
Hessian, the cost breakdown is as follows: the forward computation
calculates the gradient at cost (using the previous Section's terms)
$\le BDE$, and RM-AD on the computation of each
gradient-component would give us the Hessian at total cost $\le
B^2D^2E$. While this might not seem competitive with even the
$(D^2+D+1)E$ cost for finite-differencing, we could at least expect to
get better numerical accuracy.

If, instead, RM-AD is used for the gradient and then FM-AD for the
Hessian, the cost breakdown is as follows: the RM-AD obtains the
gradient as cost $\le BE$ and the FM-AD for the Hessian scales this by
a factor\footnote{ The effort-multiplier for computing the gradient of
a product is similar for FM-AD and RM-AD. Strictly speaking, since
there is no need to redo part of the calculation one extra time, a
better bound for the FM-AD effort multiplier is $b=B-1$ rather than
$B$. We mostly ignore this minor detail in our effort estimates.}
of~$BD$, for a total effort of~$\le B^2 D E$.  One problem with this
method is that, while it may simplify coding, the same problems
discussed in the previous Subsection are encountered for higher levels
of backpropagation, such as when computing the gradient of the lowest
eigenvalue of the Hessian (which in itself is a scalar quantity,
making RM-AD attractive).

\subsection{Option 5: `Differential Programming'}
\label{sec:diffprogramming}
It is possible to extend the ideas of the co-state equation in Eq.~\eqref{eq:ode2} to higher order terms in the spatial Taylor expansion of the loss function (i.e. higher order sensitivities).
This approach has been named `Differential Programming' in the context of machine learning model training for `Neural ODE'~\cite{chen2018neural} type
architectures in~\cite{liu2021second}.
There, the basic idea is to use an approximation for the 2nd order expansion
of the loss function to speed up optimization, much along the lines of how a dimensionally
reduced approximation to the Hessian is utilized by the L-BFGS family
of algorithms~\cite{liu1989limited}.
Both L-BFGS as well as this approach to improving `Neural ODE' training can (and
actually do) utilize heuristic guesses for the Hessian that are useful
for optimization but not entirely correct. For L-BFGS, this happens
due to using a rank-limited, incrementally-updated approximation of the
Hessian, while in~\cite{liu2021second},
practical schemes for neural ODEs with a moderate number of parameters
use both rank-limited Hessians plus also a further linearity
approximation of the form $F_i(\vec y, t)=F_{i,k}(t) y_k$
(which however is not fully explained in that article).

The `Differential Programming' approach goes as follows.
Consider a smooth scalar-valued $\mathbb{R}^D\to\mathbb{R}$ function $\mathcal{L}$
and a single-parameter-dependent diffeomorphism as
\begin{equation}
(\mathbb{R}\times \mathbb{R}^D)\to\mathbb{R}^D: \mathcal{M}=(t,\vec y)\mapsto \vec Y\,.
\end{equation}
We want to express the Taylor expansion of $\vec y\mapsto\mathcal{L}(\mathcal{M}(t, \vec y))$ in
terms of the Taylor expansion of $\vec y\mapsto\mathcal{L}(\mathcal{M}(0, \vec y))$.
In our case, we only need to work to 1st order in $t$. We have
\begin{equation}
\mathcal{M}(t, \vec y)=\vec y + t\vec F(\vec y)+\mathcal{O}(t^2)\,.
\end{equation}
The relevant partial derivatives~$M$ of~$\mathcal{M}$ at the expansion point
$(t, \vec y)=(0, \vec 0)$ are
\begin{equation} \label{eq:Mderivs}
  \begin{array}{lcl}
    M_{i, j}&=&\delta_{ij},\quad M_{i,jk\ldots}=0\,,\\
    \dot M_{i} &=& F_i(\vec 0),\quad \dot M_{i,k\ldots}= F_{i,k\ldots}(\vec 0)\,,
  \end{array}
\end{equation}
where dots indicate time-derivatives.
We introduce names for the following quantities:
\begin{equation}
\mathcal{L}(0, 0)=L\,, \qquad
\frac{\partial \mathcal{L}}{\partial t}(0, 0)=\dot L\,,\qquad
\frac{\partial \mathcal{L}}{\partial y_i}(0, 0)=L_{,i}\,.
\end{equation}

Expanding $\mathcal{L}(\mathcal{M}(t, \vec y))$ to first order in~$t$
and to second order in $y_i$ gives, schematically (where parentheses
on the right-hand-side are always for-grouping, and never evaluative):
\begin{equation}
  \begin{array}{lcccl}
    L(\mathcal{M}(t, \vec y))&=&L&+&L'(M'y)+\frac{1}{2}L''(M'y)(M'y)+\frac{1}{2}L'\underbrace{(M''yy)}_{=0}\\
    &&&+&tL'\dot M+tL'(\dot M'y)\\
    &&&+&t\Big[\frac{1}{2}L''(\dot M'y)(M'y) + \frac{1}{2}L''(M'y)(\dot M'y)\Big]\\
    &&&+&t\cdot \frac{1}{2}L'(\dot M''yy)\\
    &&&+&\{\mbox{higher order terms}\}\,,
  \end{array}
\end{equation}
where the primed notation schematically denotes partial derivatives
with respect to components of~$\vec y$ (e.g.~$L'(M' y) = L_{,i}
M_{i,j} y_j$).  In this schematic form, terms including one factor of
$t$ and no factor $y$ describe the rate-of-change for the value of the
loss $\mathcal{L}$ at the origin, terms with one factor of $t$ and one
factor $y$ describe the rate-of-change of the gradient, and-so-on. If
we consider $t$ to be an infinitesimal time-step, we can imagine doing
a suitable coordinate-adjustment after this step which brings us back
to the starting point for the next time-step again being $(t,\vec
y)=(0, \vec 0)$. Thus, we can use the above expansion to read off the
coupled evolution equations for the gradient and Hessian.

Using $M'=\delta$ and $\dot M'=F'$ as in Eq.~\eqref{eq:Mderivs},
the $[\ldots]$-term becomes
\begin{equation}
(t/2)\cdot(L_{,mj}F_{m,i}+L_{,im}F_{m,j})y_iy_j\,,
\end{equation}
and can be attributed to a tangent-space coordinate change that
will be discussed in detail later in this article.
The final term above is
\begin{equation}\label{eq:finalterm}
t \cdot \frac{1}{2}L'(\dot M''yy)=(t/2)\,L_{,m}F_{m,ij}y_iy_j
\end{equation}
which gives a contribution to the
time-evolving 2nd-order expansion coefficients (i.e. the Hessian)
coming from the gradient of the loss coupling to nonlinearities
in the velocity-function. Overall, if $\mathcal{L}$ were a nonzero
linear function (which meaningful loss functions are not), this
term arises because applying a linear function after ODE-evolution
cannot be expected to give rise to a linear function if ODE-evolution
itself is nonlinear in coordinates.
It is also interesting to note that if we start at a minimum of the loss-function and then
ODE-backpropagate the gradient and Hessian, the components of the
gradient are zero at any time $t$ and so the term~\eqref{eq:finalterm} does not
participate in the time evolution of the Hessian.
This is a nice simplification in situations where one only is interested in
backpropagating behavior around a critical point.
In this special case, it is indeed sufficient to work with the~$B_E\to B_S$ basis
transform.
In the most general situations, however, where one has a non-trivial gradient (as when using the Hessian to improve numerical optimization), the term~\eqref{eq:finalterm} does induce a contribution at $t_S$.
The contribution is proportional to (and linear in) the final-state loss-gradient, and also
proportional to the degree of nonlinearity in $F$.

Overall, using this approach, sensitivity-backpropagation of a Hessian
through an ODE works by solving the following expanded backwards-ODE:
\begin{equation}
  \label{eq:odeDDP}
\begin{split}
  (d/dt) y_k(t) &= F_k(\vec y(t))\,,\\
  (d/dt) \sigma_i(t) &= -\sigma_k(t) F_{k, i}(\vec y(t))\,,\\
  (d/dt) h_{ij}(t) &= -h_{mj}(t) F_{m, i}(\vec y(t)) -h_{im}(t) F_{m, j}(\vec y(t))\\
 & \qquad \qquad \qquad    -\sigma_m(t) F_{m,ij}(\vec y(t))\,,\\
 \end{split}
\end{equation}
where
\begin{equation}
\begin{split}
&t_0=t_E\,, \qquad t_1=t_S\,, \qquad \sigma_i(t_0)=\sigma_i(t_E)=(\partial / \partial y_i)L(\vec y(t_E))\,,\\
 &h_{ij}(t_0)=h_{ij}(t_E)=(\partial^2 / \partial y_i\partial y_j)L(\vec y(t_E)).
\end{split}
\end{equation}

So, $h_{ij}(t)$ are the coefficients of the ODE-evolving Hessian, and
$\sigma_i(t)$ are the coefficients of the gradient.
Based on the above considerations, we now pause to discuss Eq.~(9) of~\cite{liu2021second}.
There, we restrict to the case where the loss depends not on the trajectory but rather on the end-state only
$\vec y_E$ (so, $\ell=0$), and also where the ODE does not have
trainable parameters (so, $\theta=0$, meaning that all
$u$-derivatives are also zero).
Comparing with the third equation in~\eqref{eq:odeDDP} readily shows that Eq.~(9) of~\cite{liu2021second} is missing the the term as~$(d/dt) h_{ij}(t)=\ldots -\sigma_m(t) F_{m,ij}(\vec y(t))$.
This missing piece vanishes only in the following restrictive cases:
either as one considers only linear-in-$\vec y$
velocity-functions\footnote{If we obtained an autonomous system by
adding a clock-coordinate $y_D$ with $(\partial/\partial t)y_D=1$,
this corresponds to linearity in each of the first $D$ parameters, but
possibly with further, potentially nonlinear, dependency on $y_D$.},
or (as discussed above) the backpropagating starts from a critical
point (with respect to the ODE final-state coordinates) of the
scalar-valued (`loss') function~$L$ (since backpropagating
the zero gradient to any earlier time gives
$\vec\sigma(t)=\vec 0$)\footnote{The general idea that we can ignore
a contribution which in the end will merely multiply a zero gradient
has been used to great effect in other contexts, such as in particular
in~\cite{deser1976consistent} to greatly simplify the
proof~\cite{freedman1989progress} of supersymmetry
invariance of the supergravity lagrangian, where this is known
as the ``1.5th order formalism'', see~\cite{freedman2012supergravity}.}.

Intuitively, one direct way to see the need for this term is to
consider what would happen if we performed ODE-backpropagation using
the basic backpropagation recipe in which we remember all intermediate
quantities, and then backpropagating the resulting code once more: We
would encounter terms involving $F_{i,j} \equiv (\partial/\partial y_j)F_i$
on the first backpropagation, and then terms involving
$F_{i,jk}=(\partial^2/\partial y_j\partial y_k)F_i$ on the second
backpropagation.

It may be enlightening to separately discuss the role of the other
terms in the $(d/dt) h_{ij}(t)$ evolution equation, namely the contributions
\begin{equation}\label{eq:ddthterms} 
-h_{mj}(t) F_{m, i}(\vec y(t)) -h_{im}(t) F_{m, j}(\vec y(t))\,.
\end{equation}
Here, it is useful to first recall how to perform coordinate transforms of a tensorial object in multilinear
algebra\footnote{This discussion is restricted to multilinear algebra
and does \emph{not} consider changes of a position-dependent `tensor
field' with respect to general diffeomorphisms.}

Let $\vec e_i$ be the $i$-th basis vector of a (not necessarily
assumed to be orthonormal) coordinate basis $B$ and $\vec E_J$ the
$J$-th basis vector in another basis $B'$ such that these bases are
related\footnote{Einstein summation convention is understood.}
by~$\vec E_J = \vec e_i M_{iJ}$.  A vector~$\vec x=X_J \vec E_J$
in~$B'$-coordinates then has~$B$-coordinates~$\vec x=x_i\vec e_i$,
where $x_i \equiv M_{iJ}X_J$.
Given this transformation law, and
requiring that geometric relations between objects are independent of
the choice of basis, the dyadic product of two vectors $D:=\vec
v\otimes \vec w$ with coordinates $D_{IJ} = V_I W_J$ in $B'$-basis has
to have $B$-basis coordinates $d_{ij} = V_I W_J M_{iI} M_{jJ} = v_i
w_j$, so by linearity the general coordinate-transformation law is
\begin{equation} \label{eq:dtranf}
d_{ij}=D_{IJ}M_{iI}M_{jJ}\,.
\end{equation}
We get a corresponding expression for higher tensor products with one matrix-factor~$M$ per index.
Now, if~$M$ is a close-to-identity matrix, i.e.
\begin{equation} \label{eq:M}
M_{iJ}=\delta_{iJ}+\varepsilon A_{iJ}\,,
\end{equation}
we have by~\eqref{eq:dtranf} the following first-order-in-epsilon relation\footnote{We can view the coordinate transformation $M$ as
an element of the Lie group $GL(D)$, and the matrix $A$ as an element
of its Lie algebra $\mathfrak{gl}(d)$.}
\begin{equation}
  d_{ij} = D_{IJ}\delta_{iI}\delta_{jJ}+\varepsilon\left(D_{IJ}\delta_{iI}A_{jJ}+D_{IJ}A_{iI}\delta_{jJ}\right)+\mathcal{O}(\varepsilon^2).
\end{equation}
The~$\varepsilon(\cdots)$ term here directly corresponds to the structure of the terms~\eqref{eq:ddthterms}, and so Eq.~\eqref{eq:ddthterms} is interpreted as the operation transforming the coefficients of the Hessian of the objective function according to a coordinate-transformation of tangent space at a given point in time
that is induced by advancing time by $\delta t$. This is equivalent to
working out a tangent space basis $B_{E}$ at
$(t, \vec y)=(t_E, \vec y(t_E))$ and backpropagating each basis-vector
to obtain a corresponding tangent space basis $B_{S}$ at
$(t, \vec y)=(t_S, \vec y(t_S))$ which will evolve into $B_{E}$.
ODE-integrating $h_{ij}$ for a system with $F_{i,kl}=0$
then gives the same result as decomposing
the Hessian (or, using the corresponding evolution equation for higher
order terms in the Taylor expansion) with respect to $B_{E}$ and
reassembling\footnote{This is, essentially, a `Heisenberg picture'
view of the ODE.} a geometric object with the same coordinates, now
using basis $B_{S}$. Given the coordinate-linearity that results
from the $F_{i,kl}=0$ condition, it is unsurprising that these terms
codify the effect of a time-dependent $GL(D)$ transformation
uniformly acting on all of space. So, if we consider a loss $L$
as a function of $\vec y_E$ that has a convergent Taylor series
expansion, and want to obtain $L(\vec y_S)$, we
can simply basis-transform the linear and higher order terms according
to $B_E\to B_S$ and get the corresponding terms of the (still
converging) Taylor series expansion of the loss as a function of $\vec
y(t_S)$. This then implies (for example) that if we consider some $a$
such that the set $S:=\{p|L(p)\le a\}$ is convex,
the time-evolution of this set will also always be convex,
so there would be no way any such ODE could deform an ellipsoid into
a (non-convex) banana shape.

Considering the computational effort that arises when using the
complete time-evolution equation, including the $\vec\sigma F''$-term,
and also considering situations where there is no simple analytic
expression available for $F_{m,ij}$ that is easy to compute, the main
problem is that at every time-step, we have to evaluate the Jacobian
of $\sigma_m(t) F_{m,i}$, once per row $i$ of the Hessian. For this
reason, total effort here also is $\le B^2DE$, but since evolution of
the Hessian parameters is coupled, it is not straightforward to build
the Hessian row-by-row: we have to ODE-integrate a large system all in
one go.

\section{Performance comparison}

For the practitioner, it will be useful to have at least some basic
mental model for relative performance from the major methods discussed
in this work, ``differential programming'' of
Section~\ref{sec:diffprogramming} (in this Section henceforth called
``DP'') vs. the per-row backpropagation${}^{\bf (\wedge2)}$ of
Section~\ref{sec:handbackprop} (henceforth ``BP2'').

Since we clearly are limited by not being able to discuss anything
like ``a representative sample of real-world-relevant ODEs'', our
objective here is limited to answering very basic questions about
ballpark performance ratios one would consider as unsurprising, and
about what potential traps to be aware of.
Our experimental approach is characterized as follows:

\begin{itemize}
  \item For both Hessian-backpropagation methods, ODE integration is
    done in an as-similar-as-possible way, using Scientific Python
    (SciPy) ODE-integration on NumPy array data which gets forwarded
    to graph-compiled TensorFlow code for evaluating the ODE
    rate-of-change function (and its derivatives).
  \item All computations are done on-CPU, without GPU support (since
    otherwise, results may be sensitive on the degree of parallelism
    provided by the GPU, complicating their interpretation).
  \item Whenever timing measurements are done, we re-run the same
    calculation $20\times$ and take the fastest time. This in
    generally gives a better idea about the best performance
    achievable on a system than averaging.  Also, just-in-time
    compilation approaches will trigger compilation on the first
    evaluation, which should not be part of the average.
  \item The numerical problems on which timing measurements have been
    taken are low-dimensional random-parameter ODEs of the form
    \begin{equation}
    (d/dt) y_i = P^{(1)}_{ik} y_k + 0.5\,P^{(2)}_{ik\ell} y_k y_\ell
    \end{equation}
    where the coefficients $P^{(k)}$ are drawn from a normal distribution
    around 0 scaled in such a way that
    $\left\langle \left(\sum_{r_1, \ldots, r_k}P^{(k)}_{ir_1r_2\ldots r_k} x^{(1)}_{r_1} x^{(2)}_{r_2}\cdots x^{(j)}_{r_k}\right)^2\right\rangle = 1$
    for normal-distributed $x^{(j)}_{r_k}$. Calculations
    are performed in (hardware supported) IEEE-754 binary64 floating point.
    State-vector dimensionalities used were in the range $10--150$,
    in incremental steps of 10.
    Initial positions were likewise randomly-drawn from a
    per-coordinate standard normal distribution, and the loss-function was
    coordinate-distance-squared from the coordinate-origin at the
    ODE-integration endpoint.
  \item SciPy's error-tolerance management has been validated by performing
    spot-checks that compare some \codett{rtol=}\codett{atol=1e-10}
    ODE-integrations against high-accuracy results obtained via
    MPMath~\cite{johansson2013mpmath}, using 60-digit precision and
    tolerance $10^{-40}$ on a 5-dimensional problem.
  \item Per-computation, for both DP and BP2, timing measurements were
    obtained by running ODE-integrations with
    \codett{scipy.integrate.solve\_ivp()} with
    \codett{rtol=atol=1e-5} and \codett{method='RK45'}.
    Comparability of achieved result quality was 
    spot-checked at higher accuracy \codett{rtol=atol=1e-10} to agree to
    better-than-max-absolute-coefficient-distance
    $10^{-9}$. For these computations, we used
    \codett{scipy.integrate.solve\_ivp()} with \codett{method='DOP853'}.
  \item In order to validate generalizability of our results, we checked whether
    changing individual choices, such as 
    switching to SciPy's default ODE-solver `\codett{scipy.integrate.odeint()}',
    or including a 3rd-order term $\propto C_{imnp}y_my_ny_p/6$ in the ODE,
    would materially change the key insights, which was not the case.
  \item The experiments reported here have all been performed on the same
    IBM ThinkPad X1 Carbon laptop with 8 logical (according to
    \codett{/proc/cpuinfo}) Intel(R) Core(TM) i7-10610U CPU @ 1.80GHz
    cores, 16 GB of RAM, running Linux kernel version 5.19.11, with
    Debian GNU/Linux system libraries, using the public TensorFlow
    2.9.1 version with SciPy 1.8.1 and NumPy 1.21.5 (all as provided
    by PyPI, with no extra adjustments), running only default system
    processes in the background, with in total background system
    load as reported by `\codett{top}' always below 0.1 CPU-core.
\end{itemize}

On the choice of random ODEs, it should be noted that, generically,
quadratic higher-dimensional ODEs may well already exhibit the major
the known problematic behaviors of ODEs: chaos may be a generic
feature that can arise, since the three-dimensional Lorenz system can
be expressed in such a way.  Also, since $(d/dt) y(t)=-y(t)^2$ can be
embedded into these ODEs, which is solved e.g. by $y(t)=1/(t-t_0)$, it
can in principle happen that solutions ``reach infinity in finite
time''. For this reason, we kept the time-integration intervals short,
integrating from a randomly-selected starting point up to $t_{\rm
  max}=0.2$, and then backpropagating the Hessian of the loss (as a
function of the end-state) to the initial time $t=0$.

\begin{figure}
  {\includegraphics[width=0.45\textwidth]{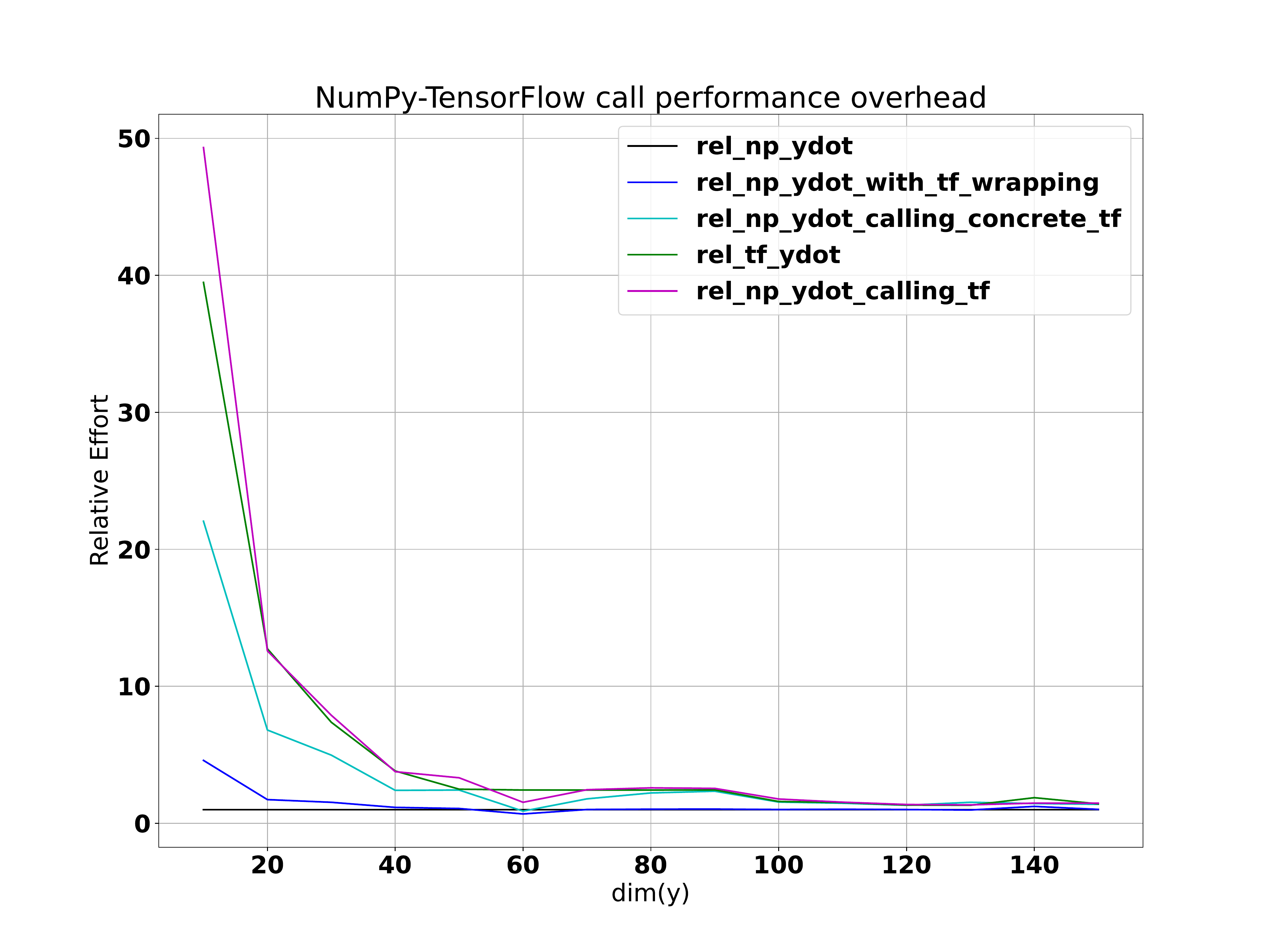}\kern1em
   \includegraphics[width=0.45\textwidth]{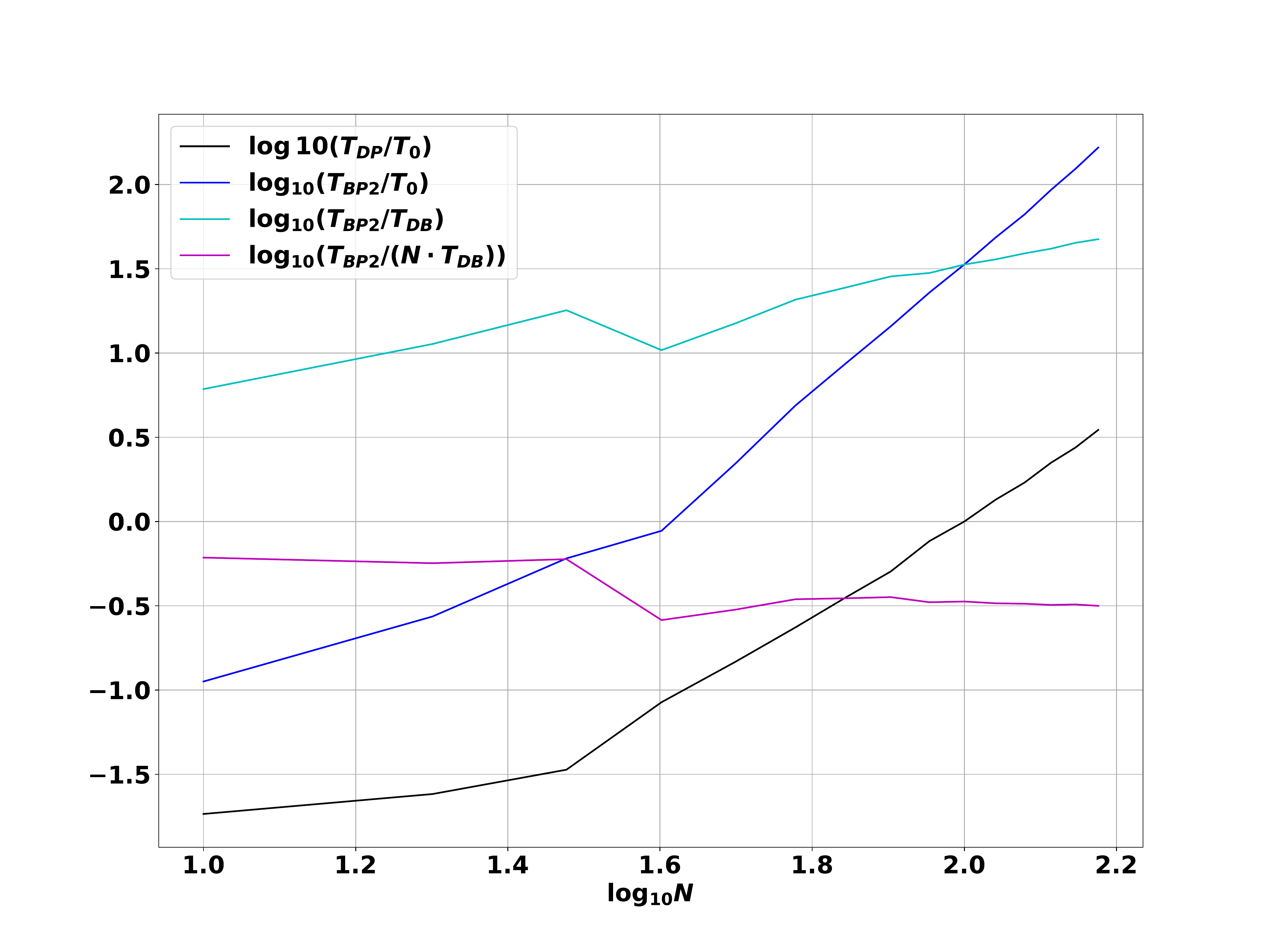}}
\begin{caption}
  {{\bf Left}: Computational latency for evaluating a function that
    computes a \codett{einsum('abc,b,c->a',\;coeffs,\;y,\;y)} tensor
    product, relative to a NumPy-only Python function.  Curve
    \codett{rel\_np\_ydot} is constant at 1,
    \codett{rel\_np\_ydot\_with\_tf\_wrapping} includes the effort to
    convert a \codett{numpy.ndarray} to a \codett{tf.Tensor} and back.
    Curve \codett{rel\_np\_ydot\_calling\_concrete\_tf} shows relative
    effort for calling the `concrete TensorFlow function' for the
    given input-signature. Curve \codett{rel\_tf\_ydot} shows time for
    TensorFlow to compute a \codett{tf.Tensor} result from a constant
    \codett{tf.Tensor} input. Curve
    \codett{rel\_np\_ydot\_calling\_tf} shows relative effort for
    calling a just-in-time-compiled \codett{@tf.function}-decorated
    function \codett{tf\_f} not-first-time via
    \codett{tf\_f(tf.constant(data,\;dtype=tf.float64)).numpy()}.
    This in particular includes the effort to determine (through
    Python-level code) the specific compiled realization from the
    call-signature. Clearly, for small vector-sizes, the overhead can
    be very sizeable.  {\bf Right}:
    $\log_{10}(\mbox{effort-ratio})$/$\log_{10}(\mbox{time / reference\_time})$
    log/log-plot for computing the Hessian via
    methods `BP2' and `DP'.  Reference time $T_0$ is computation of a
    Hessian via method `DP' for state vector size 100: The point
    $(\log_{10} N, 0)$ sits on the black curve.  As the slope of the
    black curve indicates, computational effort grows roughly
    cubically for problem sizes above $N=40$. The `BP2' curve sits
    about 1.5 units above the `DP' curve, indicating a typical
    relative-effort factor of about 30 in the range shown
    here. Overall, while the raw effort-ratio (cyan) grows with
    increasing problem size, growth is sublinear, as is also indicated
    when dividing the effort-ratio by the problem size (cyan). This
    suggests that the effort-ratio should (as expected) approach a
    constant for much larger problem sizes, when the overhead to
    forward numerical data between the ODE-integrator and the RM-AD
    framework is fully amortized.}
\label{fig:tf_overhead}
\end{caption}
\end{figure}

One important insight is that, while one would on theoretical grounds
expect equivalent scaling with problem size, the `BP2' approach suffers
from having to do many independent ODE-integrations, one per
state-space dimension. This means that if the coordinate-derivatives
of the rate-of-change function have been obtained from some automatic
differentiation framework (such as in the case at and TensorFlow), any
overhead in exchanging numerical data between this framework and the
ODE integration code (such as: constant-effort-per-call bookkeeping)
easily can become painful, and ultimately would amortize only at
comparatively large problem-size. This is shown in
figure~\ref{fig:tf_overhead}-left. Henceforth, timing measurements have
been performed with TensorFlow concrete-function resolution being
taken out of the ODE integrator's main loop.

Figure~\ref{fig:tf_overhead}-right shows relative performance in a
log/log-plot, where slopes correspond to scaling exponents\footnote{As
a reminder, if $y=cx^\alpha$, we have $(\log y)=\alpha(\log x)+\log
c$, so power-law scaling behavior shows as a straight line in a
log/log plot.}.

For this problem, where effort~$E$ is observed to scale roughly
as~$E\propto N^3$, so like to the number of multiplications for
evaluating the right hand side, we find that method `DP2' is at a
major disadvantage due to having to perform many individual
ODE-integrations which make any constant-effort overhead for RM-AD to
ODE-integrator adapter code painful. Still, it appears plausible that
for very large problem sizes, the effort-ratio approaches a constant,
as suggested by theoretical considerations.

Overall, we conclude that despite method `BP2' being expected to have
asymptotically-same time-complexity (i.e. scaling behavior with
problem size) as method `DP' if we consider equally-sized ODE
time-steps for fixed step size integrators, its big achilles heel is
that a large number of (small) ODE integrations that to be performed,
and this can prominently surface any relevant constant-effort function
call overhead of evaluating the rate-of-change function (and its
derivatives). For not-implausible tool combinations and smallish ODEs,
observing a performance-advantage of `DP' over `BP2' by a factor 30
should not come as a surprise. This however may look rather different
if highly optimized code can be synthesized that in particular avoids
all unnecessary processing.

\section{Hand-crafting an ``orbit non-closure'' Hessian}
\label{sec:handcrafted}

We want to illustrate some of the subtleties that can arise when
hand-backpropagating ODEs to compute Hessians by studying an
`orbit-nonclosure' loss. This example illustrates how to address
subtleties that arise when an objective function depends on both
initial and final state of the ODE-integration.

The idea is that for some dynamical system,
such as a collection of planets, the state-of-motion is completely
described by a collection of positions and velocities (equivalently,
momenta), one per object. We are then asking, for a given time
interval $\Delta T$, how large the misalignment is between the
state-of-motion vector after ODE-integrating the equations of motion
over time-interval $\Delta T$ and the starting
state-of-motion. Somewhat unusually, the loss here depends on
\emph{both} the initial- and final-state of ODE-integration, in a
non-additive way. Overall, this detail complicates
hand-backpropagation, and we will address the need to make the
construction's complexity manageable by introducing an appropriate
formalism. This, then, is not tied to
backpropagation-of-backpropagation applications, but should be
generically useful for reasoning about the `wiring' of complicated
field theoretic constructions (including `differential programming').

The code that is published alongside this article~\footnote{It is
available both as part of the arXiv preprint and also on:
\url{https://github.com/google-research/google-research/tree/master/m\_theory/m\_theory\_lib/ode}}
provides a Python reference implementation for the algorithms
described in Sections~\ref{sec:handbackprop}
and~\ref{sec:diffprogramming} in general terms in a form that should
admit straightforward implementation in other software frameworks in
which there is only limited support for (ODE-)backpropagation.

Using the grammar~described in Appendix~\ref{app:tcd_dsl}, we are given
the function \codett{ydot} that computes the ``configuration-space
velocity'' as a function of configuration-space position. We explore
this for autonomous systems, i.e. \codett{ydot} is time-independent
and has signature:
\begin{texteq} \label{eq:ydotsignature}
\begin{verbatim}
    ydot: ^(y[:]) => [:]
\end{verbatim}
\end{texteq}
The ODE-integrator \codett{ODE} computes a final configuration-space
position given the velocity-function \codett{ydot}, the initial
configuration-vector \codett{y[:]}, as well as the starting time
\codett{t0[]} and final time \codett{t1[]} and hence has signature:
\begin{texteq*}
\begin{verbatim}
    ODE: ^(ydot{}, y[:], t0[], t1[]) => [:]
\end{verbatim}
\end{texteq*}
The ODE-backpropagation transformation \codett{OBP} turns a function
with the signature of \codett{ydot}, as in Eq.~\eqref{eq:ydotsignature}, into a velocity-function for the
backpropagating ODE-integration. There are multiple ways to write this
given the grammar, but the idea likely is most clearly conveyed by
this form, where the position-vector slot-name for the `doubled' ODE
is also \codett{y[:]} (to allow forming the doubled-ODE of a doubled-ODE),
and the configuration-space dimensionality is \codett{dim}.
\begin{texteq} \label{eq:OBP}
\begin{verbatim}
    OBP = ^ (ydot{},
             y[:],
             vy[:] = ydot(y[:]=y[:dim]),
             J_ydot[:, :] = ydot(y[:]=y[:dim])[:, d y[:]]) ->
           &concat(0, vy[:], &es(-J_ydot[:, :] @ i, j; y[dim:] @ i -> j))
\end{verbatim}
\end{texteq} 
According to the language rules in appendix~\ref{app:tcd_dsl}, we can
alternatively regard \codett{OBP} as a higher order function that maps
a function-parameter named `\codett{ydot}` to a function of a
configuration space position `\codett{y[:]}`, hence another function
with the signature of `\codett{ydot}`, but for the twice-as-wide
vector, incorporating also the costate equation.

For the example problem of orbit non-closure, we need a `target'
function that measures the deviation between configuration space
positions. We will subsequently evaluate this on the initial and final
state of ODE-integration for a given time interval and initial
position, but the most straightforward ``L2 loss'' target function
one can use here is:
\begin{texteq} \label{eq:Tloss}
\begin{verbatim}
    T = ^(pos0[:], pos1[:],
          delta[:] = pos0[:] - pos1[:]) ->
        &es(delta[:] @ a; delta[:] @ a ->)
\end{verbatim}
\end{texteq} \
It makes sense to introduce names for the first and second derivatives of \codett{T}:
\begin{texteq} \label{eq:Tpartials}
\begin{verbatim}
    T0 = ^(pos0[:], pos1[:]) -> T(pos0[:]=pos0[:],pos1[:]=pos1[:])[d pos0[:]]
    T1 = ^(pos0[:], pos1[:]) -> T(pos0[:]=pos0[:], pos1[:]=pos1[:])[d pos1[:]]
    T00 = ^(pos0[:], pos1[:]) -> T0(pos0[:]=pos0[:], pos1[:]=pos1[:])[:, d pos0[:]]
    T01 = ^(pos0[:], pos1[:]) -> T0(pos0[:]=pos0[:], pos1[:]=pos1[:])[:, d pos1[:]]
    T11 = ^(pos0[:], pos1[:]) -> T1(pos0[:]=pos0[:], pos1[:]=pos1[:])[:, d pos1[:]]
\end{verbatim}
\end{texteq}
With these ingredients, we can formulate orbit non-closure \codett{NC}
as a function of the velocity-function, initial configuration-vector,
and orbit-time as follows:
\begin{texteq} \label{eq:NC}
\begin{verbatim}
    NC = ^(ydot{}, y_start[:], t_final[],
           y_final[:] = ODE(ydot{}=ydot, y[:]=y_start[:],
                            t0[]=0, t1[]=t_final[])[:]) ->
         T(pos0[:]=y_start[:], pos1[:]=y_final[:])[]
\end{verbatim}
\end{texteq} 
The gradient of this orbit non-closure with respect to the initial position,
which we would write as:
\begin{texteq} \label{eq:gradNC}
\begin{verbatim}
    grad_NC = ^(ydot{}, y_start[:], t_final[]) ->
              NC(ydot{}=ydot, y_start[:]=y_start[:],
                 t_final[]=t_final[])[d y_start[:]]
\end{verbatim}
\end{texteq} 
clearly is useful for finding closed orbits via numerical minimization, 
and in case the original ODE is well-behaved under time-reversal,
this can be computed as follows:
\begin{texteq} \label{eq:gradNCbackprop}
\begin{verbatim}
    grad_NC = ^(
      ydot{}, y_start[:], t_final[],
      y_final[:] = ODE(ydot{}=ydot, y[:]=y_start[:], t0[]=0, t1[]=t_final[])[:],
      s_T_y_start[:] = T(pos0[:]=y_start[:], pos1[:]=y_final[:])[d pos0[:]],
      s_T_y_final[:] = T(pos0[:]=y_start[:], pos1[:]=y_final[:])[d pos1[:]],
      s_T_y_start_via_y_final[:] = ODE(ydot{}=OBP(ydot{}=ydot),
                                       y[:]=&concat(0, y_final[:], s_T_y_final[:]),
                                       t0[]=t_final[], t1[]=0)[dim:]) ->
      s_T_y_start[:] + s_T_y_start_via_y_final[:]
\end{verbatim}
\end{texteq} 
The task at hand is now about computing the gradient (per-component, so Jacobian)
of \codett{grad\_NC} with respect to \codett{y\_start[:]},
 \begin{texteq} \label{eq:hessianNC}
\begin{verbatim}
    hessian_NC = ^(ydot{}, y_start[:], t_final[]) ->
                  NC(ydot{}=ydot, y_start[:]=y_start[:],
                     t_final[]=t_final[])[d y_start[:], d y_start[:]]
\end{verbatim}
\end{texteq} 
which is the same as
 \begin{texteq} \label{eq:hessianNCbackprop}
\begin{verbatim}
    hessian_NC = ^(ydot{}, y_start[:], t_final[]) ->
                  grad_NC(ydot{}=ydot, y_start[:]=y_start[:],
                           t_final[]=t_final[])[:, d y_start[:]]
\end{verbatim}
\end{texteq} 

Let us consider a fixed component of~\codett{grad\_NC}.
Let~\codett{j} henceforth be the index of this component.
We first proceed by introducing the following auxiliary vector-valued function~\codett{S\_fn},
which computes the sensitivity of the~\codett{grad\_NC[j]} component
with respect to the starting position of the ODE-integration for~\codett{s\_T\_y\_start\_via\_y\_final[:]} in Eq.~\eqref{eq:gradNCbackprop}.
We can now write this as:
 \begin{texteq} \label{eq:Sn}
\begin{verbatim}
    S_fn = ^( 
      ydot{}, y_start[:], t_final[],
      y_final[:] = ODE(ydot{}=ydot, y[:]=y_start[:], 
                       t0[]=0, t1[]=t_final[])[:],
      s_T_y_final[:] = T1(pos0[:]=y_start[:], pos1[:]=y_final[:])[:]) ->
         ODE(ydot{}=OBP(ydot{}=ydot),
             y[:]=&concat(0, y_final[:], s_T_y_final[:]),
             t0[]=t_final[], t1[]=0)[dim:][j, d y[:]]
\end{verbatim}
\end{texteq}

Using ODE backpropagation again to obtain the gradient, and remembering
that we can regard taking the $j$-th vector component as a scalar product
with a one-hot vector, i.e.

\begin{texteq*}
\begin{verbatim}
    x[j] = &es(x[:] @ k; &onehot(j, dim) @ k ->)
\end{verbatim}
 \end{texteq*}
we can compute~\eqref{eq:Sn} as:
\begin{texteq}
\begin{verbatim}
    S_fn = ^ ( 
      ydot{}, y_start[:], t_final[],
      y_final[:] = ODE(ydot{}=ydot, y[:]=y_start[:], t0[]=0, t1[]=t_final[])[:],
      s_T_y_final[:] = T1(pos0[:]=y_start[:], pos1[:]=y_final[:])[:],
      ydot2{} = OBP(ydot{}=ydot),
      obp1_start[:] = &concat(0, y_final[:], s_T_y_final[:])[:],
      s_T_y_start_via_y_final_ext[:] = ODE(ydot{}=ydot2,
                                           y[:]=obp1_start[:],
                                           t0[]=t_final[], t1[]=0)[:] ) ->
        ODE(ydot{}=OBP(ydot{}=ydot2),
            t0[]=t_final[], t1[]=0,
            y[:]=&concat(0,
                         s_T_y_start_via_y_final_ext[:],
                         &onehot(dim + j, 2 * dim)))
\end{verbatim}
\end{texteq}

Here, unlike in the initial version of \codett{S\_fn}, we have not
trimmed \codett{s\_T\_y\_final} to its second half\footnote{Since this
is an `internal' definition, both forms still describe the same field
theoretic function}, to show more clearly the conceptual form of this
2nd backpropagation approach.

At this point, we encounter a relevant detail: the second
\codett{ODE(...)}  starts from an initial position that contains the
reconstructed starting-position of the system \codett{y\_start[:]} obtained
after numerically going once forward and once backward through
ODE-integration. This is exactly the code that an algorithm which
provides backpropagation as a code transformation would synthesize
here. Pondering the structure of the problem, we can however eliminate
some numerical noise by instead starting from the known good starting
position, i.e. performing the second ODE-backpropagation as
follows\footnote{On the toy problems used in this work to do performance
measurements, this replacement does however not make a noticeable difference.}:
\begin{texteq}
\begin{verbatim}
    ODE(ydot{}=OBP(ydot{}=ydot2),
        t0[]=t_final[], t1[]=0,
        y[:]=&concat(0,
                     y_start[:],
                     s_T_y_start_via_y_final_ext[dim:],
                     &zeros(dim),
                     &onehot(j, dim)))
\end{verbatim}
\end{texteq}

Now, with \codett{S\_fn}, we have the sensitivity of an entry of the
orbit-nonclosure gradient on the initial state of the
ODE-backpropagation in Eq.~\eqref{eq:gradNCbackprop}.
At this point, it might be somewhat difficult to
see what to do when not having a formalism that helps reasoning this
out, but going back to Eq.~\ref{eq:gradNC}, it is clear where we are:
we need to find the sensitivity of \codett{gradNC[j]} on the input
that goes into \codett{T0} and \codett{T1}.

Once we have these sensitivities of the $j$-th entry of the gradient
on \codett{y\_start[:]} and \codett{y\_final[:]} via the
\codett{pos0[:]} respectively \codett{pos1[:]} slots of \codett{T0},
\codett{T1}, we can backpropagate the sensitivity on
\codett{y\_final[:]} into an extra contribution to the full
sensitivity on \codett{y\_start[:]}.

With this insight, we are ready to complete the construction,
obtaining the $j$-th row of the Hessian now as follows.
We use the prefix \codett{s\_} as before to indicate the sensitivity of the value
of the target function \codett{T[]} on various intermediate
quantities, and the prefix \codett{z\_} to indicate the sensitivity of
\codett{gradNC[j]} on some specific intermediate quantity.
Comments~\codett{F\{nn\}} label the forward-computation steps, and on the
backward pass label reference these forward-labels to flag up the use
of an intermediate quantity that produces a contribution to the
sensitivity.
\begin{texteq} \label{eq:hessianrowj}
\begin{verbatim}
    hessian_row_j = ^(
      ydot{}, y_start[:], t_final[],  # F1
      y_final[:] = ODE(ydot{}=ydot, # F2
                       y[:]=y_start[:],
                       t0[]=0,
                       t1[]=t_final[])[:],
      # Value for which we obtain the Hessian - not needed here:
      # T_val = T(pos0[:]=y_start[:], pos1[:]=y_end[:])
      s_T_y_start[:] = T0(pos0[:]=y_start[:], pos1[:]=y_final[:]),  # F3
      s_T_y_final[:] = T1(pos0[:]=y_start[:], pos1[:]=y_final[:]),  # F4
      obp1_start[:] = &concat(0, y_final[:], s_T_y_final[:])[:],  # F5
      ydot2{} = OBP(ydot{}=ydot),
      s_T_y_start_via_y_final[:] = ODE(ydot{}=ydot2, # F6
                                       y[:]=obp1_start,
                                       t0[]=t_final[], t1[]=0)[dim:],
      # The j-th entry of the gradient of T_val.
      # Not needed, but shown to document structural dependency.
      # grad_T_val_entry_j[] = (
      #  [s_T_y_start[:] + s_T_y_start_via_y_final[:]][j]),  # F7
      onehot_j[:] = &onehot(j, dim),
      # We want to know the sensitivity of grad_T_val_entry_j[] on
      # y0. Subsequently, let zj_{X} be the sensitivity of this quantity
      # on the corresponding intermediate quantity {X}.
      # Processing dependencies of intermediate quantities in reverse,
      # obp1_start, then s_T_y_final, s_T_y_start, finally y_final:
      z_obp1_start = ODE(ydot{}=OBP(ydot{}=ydot2),
                         y[:]=&concat(0,
                                      y_start[:],
                                      s_T_y_start_via_y_final[:],
                                      &zeros(dim),
                                      onehot_j))[2*dim:],
      z_s_T_y_final[:] = z_obp1_start[dim:],  # from F5
      z_s_T_y_start[:] = onehot_j[:],  # from F7
      z_y_final[:] = (
        # from F5
        z_obp1_start[:dim] +
        # from F4
        &es(T11(pos0[:]=y_start[:], pos1[:]=y_final[:]) @ a, b;
            z_s_T_y_final[:] @ b -> a) +
        # from F3
        &es(T01(pos0[:]=y_start[:], pos1[:]=y_final[:]) @ a, b;
            z_s_T_y_start[:] @ b -> a))
      ) ->
      # The result is "z_y_start[:]", i.e. grad_T_val[j, d y_start[:]].
      ( # from F3
        &es(T00(pos0[:]=y_start[:], pos1[:]=y_final[:]) @ a, b;
            z_s_T_y_start @ a -> b) +
        # from F4
        &es(T01(pos0[:]=y_start[:], pos1[:]=y_final[:]) @ a, b;
             z_s_T_y_start @ a -> b) +
        # from F2
        ODE(ydot{}=ydot2,
            y[:]=&concat(0, y_final[:], z_y_final[:])[:],
            t0[]=t_final[],
            t1[]=0)[dim:])
\end{verbatim}
\end{texteq}

We can then readily stack these Hessian-rows into the
Hessian~$h_{jk}$. Overall, while the Hessian always is a symmetric
matrix, this numerical approach gives us a result that is not
symmetric-by-construction, and numerical inaccuracies might actually
make $h_{jk}\neq h_{kj}$. The magnitude of asymmetric part can be
regarded as providing an indication of the accuracy of the result, and
one would in general want to project to the symmetric part.

\section{Hessians for ODEs -- examples}
\label{sec:HessianODEExamples}

In this Section we apply the backpropagation of gradients and
Hessians, as developed in Section~\ref{sec:handcrafted}, through ODEs
describing specific dynamical systems of interest.  We focus on
Hamiltonian systems, where the time-evolution for each component of
phase-space is given by a first-order ODE\footnote{For such
Hamiltonian systems, one would in general want to use not a generic
numerical ODE-integrator, but a symplectic integrator. We here do not
make use of such special properties of some of the relevant ODEs,
since the achievable gain in integrator performance does not play a
relevant role here.} in wall-time.  More precisely, we consider
Hamiltonian functions~\texttt{H\{\}} with signature
\begin{texteq}
\begin{verbatim}
    H: ^(y[:]) => []
\end{verbatim}
\end{texteq}
where~\texttt{y[:]} is a~\texttt{d+d}-dimensional phase space vector, which should be treated as the concatenation of the position~\texttt{d}-vector and the conjugate momentum~\texttt{d}-vector:
\begin{texteq} \label{phasespacevector}
\begin{verbatim}
    y[:] = &concat(0, q[:], p[:]) 
\end{verbatim}
\end{texteq}
Given such Hamiltonian function, the time-evolution ODE of~\texttt{y[:]} is given by Hamilton's equations
as\footnote{For many physical systems, these amount to: position rate-of-change, i.e. velocity,
is proportional to momentum, and momentum rate-of-change is proportional to the spatial gradient
of total energy, i.e. the change in potential energy.}
\begin{texteq} \label{HamiltonsEqn}
\begin{verbatim}
    ydot_H = ^(y[:], q[:]=y[:d], p[:]=y[d:]) 
               -> &concat(0, H[d q[:]], -H[d p[:]])[:]
\end{verbatim}
\end{texteq}

In this Section, the Hamiltonian function~\texttt{H\{\}} will successively be taken to be that of the 3d isotropic harmonic oscillator, the Kepler problem, and the 2d three-body-problem.
In each case, the loss function~\texttt{NC\{\}} is defined in Eq.~\eqref{eq:NC}, Eq.~\eqref{eq:Tloss}, to be the squared-coordinate-distance misalignment between an initial configuration \texttt{y\_start[:]} at an initial time (that we fix to \texttt{0}) and the corresponding time-evolved end-state configuration \texttt{y\_final[:]} at boundary time \texttt{t1[]}.
For convenience, we copy it here:
\begin{texteq} \label{Loss}
\begin{verbatim}
    NC = ^(ydot{}, y_start[:], t_final[],
           y_final[:] = ODE(ydot{}=ydot, y[:]=y_start[:],
                            t0[]=0, t1[]=t_final[])[:])
                     ->
                     T(pos0[:]=y_start[:], pos2[:]=y_final[:])[]
\end{verbatim}
\end{texteq}
where
\begin{texteq} \label{T_fn_examples}
\begin{verbatim}
    T = ^(pos0[:], pos1[:],
          delta[:] = pos0[:] - pos1[:]) ->
        &es(delta[:] @ a; delta[:] @ a ->)
\end{verbatim}
\end{texteq}

For a given Hamiltonian system as in Eq.~\eqref{HamiltonsEqn}, the initial configurations of interest on these systems are configurations \texttt{y\_start[:]=Y0[:]} with fixed boundary time~\texttt{t1[]=T\_final[]} corresponding to  closed orbits.
For a given~\texttt{t1[]=T\_final[]}, we refer to~\texttt{Y0[:]} as a~\textit{solution}.
Our exact task here is then to perform (using our backpropagation formalism) the evaluation of the loss function, its gradient function and Hessian functions on such solutions.
That is, we wish to compute:
\begin{texteq} \label{functionCall}
\begin{verbatim}
    # 0-index loss
    loss = NC(ydot{}=ydot_H, y_start[:]=Y0[:],
              t_final[]=T_final)
    # 1-index gradient
    grad = grad_NC(ydot{}=ydot_H, y_start[:]=Y0[:],
                   t_final[]=T_final)
    # 2-index Hessian
    hessian = hessian_NC(ydot{}=ydot_H, y_start[:]=Y0[:],
                         t_final[]=T_final)
\end{verbatim}
\end{texteq}
If~\texttt{y\_start[:]=Y0[:]} with~\texttt{t1[]=T\_final[]} is indeed a solution,
then the loss should be zero. If the solution is stationary, we
also expect the gradient-vector-function to evaluate to the zero-vector.
Finally, by analyzing the flat directions of the Hessian, we comment on the symmetries of the solution.
These are all independent deformations  of an orbit-closing solution \texttt{Y0[:] -> Y0\_def[:]},
such that \texttt{Y0\_def[:]} is also an orbit-closing solution for the prescribed orbit-time.

Continuous-parameter symmetries that act on phase-space in a way that
commutes with time-translation, i.e. the action of the Hamiltonian, will
neither affect time-to-orbit nor orbit-closure. Hence, their infinitesimal
action-at-a-point belongs to the nullspace of the Hessian of the
aforementioned function, evaluated at that point. It can however happen
that a Hessian's nullspace is larger than deformations that correspond to
symmetries of the system. This can happen for different reasons: there can
be deformations that do not affect orbit-closure and also time-to-orbit,
but nevertheless change the total energy, i.e. the value of the
Hamiltonian, or the higher (e.g. fourth) order corrections from the Taylor
expansion of the loss function might not respect the symmetry of the
system.

Hence, the dynamical symmetries of the Hamiltonian always form a (not
necessarily proper) subspace of the symmetries of a solution that leave
orbit-time\footnote{For example,~in our context of fixed orbit-time, the
``deformation'' on~\texttt{Y0[:]} generated by the time-translational
symmetry of the Hamiltonian is simply the act of shifting
the~\texttt{Y0[:]} along the orbit.}
invariant\footnote{An interesting exercise would be tweaking the
misalignment function~\texttt{T\{\}} in~Eq.~(\ref{T_fn_examples}) in way
that does not allow for such dynamical symmetry-violating deformations. One
would then want to use coordinate-misalignment loss to find a closed orbit,
and an extended loss function $\tilde L(\vec y_{\rm start},\vec y_{\rm
  end})=L(\vec y_{\rm start},\vec y_{\rm end})+(E_{\rm start}-E_0)^2$ which
also punishes energy-changing deformations to find candidates for
generators of dynamical symmetries.}.

\subsection{The 3d Harmonic Oscillator}

As a first simple example, consider the classical isotropic 3d harmonic oscillator. 

The harmonic oscillator describes the dynamics of a point mass \texttt{m[]} in
three-dimensional space, experiencing a force towards the coordinate-origin that
is proportional to its distance from that point, with spring constant \texttt{k[]}.
A generic state-vector \texttt{y[:]} is the concatenation of the
six phase space coordinates contained in \texttt{q[:]} and \texttt{p[:]}, which respectively
correspond to position and conjugate momentum vectors of the mass.
For simplicity we set the mass and spring constant to unity.
The Hamiltonian function for this system is then
\begin{texteq} \label{hoscHamil}
\begin{verbatim}
    H_ho3 = ^ (y[:]) ->  (1/2) * &es(y[:3] @ qi; y[:3] @ qi ->)
                      +  (1/2) * &es(y[3:] @ pi; y[3:] @ pi ->)
\end{verbatim}
\end{texteq}
The rate-of-change function for a phase-space vector \texttt{y[:]} in
this system is given by Hamilton's equations~\eqref{HamiltonsEqn} as
\begin{texteq} \label{hoscsystem}
\begin{verbatim}
    ydot_ho3 = ^(y[:]) -> &concat(0, y[3:], -y[:3])[:]
\end{verbatim}
\end{texteq}

We now come to the task of evaluating for the harmonic oscillator Eq.~\eqref{hoscHamil} the loss function Eq.~\eqref{Loss}, its gradient function and its Hessian functions.
As discussed in the introduction of this Section, we choose a particular initial configuration \texttt{y\_start[:]=Y0[:]} and fixed boundary
time \texttt{t\_final[]=T\_final[]} such that the corresponding trajectory is a closed orbit. 
In the harmonic oscillator, such configurations are simply obtained by setting \texttt{T\_final[]} to be the ``orbit time'', which is a function of the mass and spring constant only, $T=2\pi\sqrt{m/k}$.
For this choice of boundary time, any initial state vector yields an orbit.
For \texttt{m[] = k[] = 1}, we hence have:
\begin{texteq}  \label{hoscTfinal}
\begin{verbatim}
    T_final[] = 6.28318530718  # 2 pi
\end{verbatim}
\end{texteq}
The arbitrary initial configuration \texttt{Y0[:]} is chosen as
\begin{texteq} \label{hoscSolution}
\begin{verbatim}
    Y0[:] = [50., 10., 50., -20., 10., -0.1]
\end{verbatim}
\end{texteq}
We now have all the data to perform an evaluation the loss, gradient and Hessian functions.
The loss and its gradient for \texttt{y\_start[:]=Y0[:]} and~\texttt{t1[]=T\_final[]} given in Eq.~\eqref{hoscTfinal} are indeed numerically-zero:
\begin{texteq} \label{hoscOutputLossGrad}
\begin{verbatim}
    NC(ydot{}=ydot_ho3, y_start[:]=Y0[:], t_final[]=T_final[])[]
    # = 1.0523667647935759e-17
                                  
    grad_NC(ydot{}=ydot_ho3, y_start[:]=Y0[:],
            t_final[]=T_final[])[:] 
    #    = [1.22409460e-09, -9.60937996e-10, -4.89819740e-10,
    #       4.50560833e-09, 7.61573915e-10,  4.30725394e-09]
\end{verbatim}
\end{texteq}
The vector of eigenvalues of the Hessian also has numerically-zero entries:
\begin{texteq}
\begin{verbatim}
    &eigvals(hessian_NC(ydot{}=ydot_ho3, y_start[:]=Y0[:],
             t_final[]=T_final[])[:,:])[:]
    #   = [-5.9e-11, -5.9e-11, -5.9e-11,
    #      4.4e-11, 4.4e-11, 4.4e-11]
\end{verbatim}
\end{texteq}
The zero loss and gradient in Eq.~\eqref{hoscOutputLossGrad} indicates that
the initial configuration Eq.~\eqref{hoscSolution} with boundary time
Eq.~\eqref{hoscTfinal} is indeed a stable orbit-closure solution for the
harmonic oscillator.  The zero Hessian eigenvalues (in fact, all entries of
the Hessian are zero) captures the fact that with boundary time
Eq.~\eqref{hoscTfinal}, any arbitrary deformation of~\texttt{Y0[:]}
preserves orbit closure.
It is therefore interesting to note that while the dynamical symmetry group
of the Hamiltonian Eq.~\eqref{hoscHamil} is~$U(3)$, the symmetries for
orbit closure are in actually~$GL(6)$ (i.e. the group of all invertible
transformations on the 6-vector~\texttt{Y0[:]}). This is then an example of
the aforementioned phenomenon that the orbit closure criterion can admit
additional symmetries beyond those of the dynamical symmetries of the
system. The transformations on phase space associated to these extra
symmetries will excite the system to different values for the conserved
charges of the Hamiltonian. We will see shortly that, for instance, they
may map one orbit-closure solution with energy
$E=\texttt{H\_ho3(y[:]=Y0[:])}$ to another orbit-closure solution with
energy $E' \neq E$.

In Figure~\ref{fig:orbitplots}(a), the full line shows the orbit whose initial configuration is the solution~\texttt{Y0[:]}.
The dotted line shows the trajectory whose initial configuration is a small (and arbitrarily chosen) deformation of~\texttt{Y0[:]}:
\begin{texteq} \label{hoscSolutionDef}
\begin{verbatim}
    Y0_deformed[:] =  [50.70710678 , 10. ,  50. , -20. , 9.29289322,  -0.1]
\end{verbatim}
\end{texteq}
This trajectory is also an orbit with the same orbit time Eq.~\eqref{hoscTfinal}, 
as expected from the discussion in the prior paragraph.
The loss on this deformed solution is
\begin{texteq}
\begin{verbatim}
    NC(ydot{}=ydot_ho3, y_start[:]=Y0_deformed[:], t_final[]=T_final[])[]
    #  = 1.0717039249337393e-17
\end{verbatim}
\end{texteq}
Finally, note that the original solution~\texttt{Y0[:]}
and the deformed solution~\texttt{Y0\_deformed[:]} have differing energies,
explicitly illustrating the aforementioned situation:
\begin{texteq}
\begin{verbatim}
    H_ho3(y[:]=Y0[:]) # = 2800.
    H_ho3(y[:]=Y0_deformed[:]) # = 2828.79
\end{verbatim}
\end{texteq}

\subsection{The 3d Kepler Problem}

A more complex system is the 3d Kepler problem. This is the gravitational
two-body problem where the central force between the two masses is taken to
be proportional to the squared-inverse Euclidean separation.  This system
is naturally symmetric under three-dimensional spatial translations, and we
use this to fix the origin in the centre of mass of the two bodies.
This procedure reduces the problem to the motion of one body, with mass
equal to the so-called reduced mass~\texttt{mu[] = m1[]*m2[]/(m1[]+m2[])},
in a central field with potential energy proportional to the inverse
distance of the body to the origin.
For simplicity, we also choose the masses \texttt{m1[]}, \texttt{m2[]} such
that~\texttt{mu[]=1}. The Hamiltonian function in terms of a generic
phase-space vector
\begin{texteq*}
\begin{verbatim}
    y[:] = &concat(0, q[:], p[:])[:]
\end{verbatim}
\end{texteq*}
is given, for an appropriately chosen Newton's constant, as
\begin{texteq} \label{keplHamil}
\begin{verbatim}
    H_kepl = ^ (y[:]) -> (1/2) * &es(y[3:] @ i; y[3:] @ i ->)
                               - &pow(&es(y[:3] @ i; y[:3] @ i ->), -1/2)
\end{verbatim}
\end{texteq}
The rate-of-change vector-function of \texttt{y[:]} is given by Hamilton's equations as
\begin{texteq} \label{keplsystem}
\begin{verbatim}
    ydot_kepl = ^(y[:],
                  r_factor[] = &pow(&es(y[:3] @ i; y[:3] @ i ->), -3/2))
                  -> &concat(0, y[3:], -y[:3] * r_factor[])
\end{verbatim} 
\end{texteq}
As for the prior example (harmonic oscillator), we evaluate the
orbit-closure loss function given in Eq.~\eqref{Loss}, as well as its
gradient and Hessian functions, on an orbit-closure
configuration~\texttt{(y\_start[:]=Y0[:], t1[]=T\_final[])}.
Unlike for the harmonic oscillator, we do not use analytical expressions
of the system to ab-initio prescribe a good \texttt{T\_final[]} and \texttt{Y0[:]}.
Instead, we first fix \texttt{T\_final[]} to a sensible value and then let a BFGS
optimizer (that minimizes the orbit-closure loss~\texttt{NC\{\}}) find a solution for \texttt{Y0[:]}
based on some initialization \texttt{Y0\_init[:]}. Note that the BFGS-type
optimizers rely on knowing a gradient at each iteration, which we provide by
allowing the optimizer to evaluate the \texttt{grad\_NC\{\}} function.
Note also that throughout this Section, our optimizer runs are set with target
gradient norm tolerance of \texttt{1e-12}. Now, for the boundary time, we fix
\begin{texteq} \label{keplTfinal}
\begin{verbatim}
    T_final[] = 6.28318530718
\end{verbatim}
\end{texteq}
while for the state-vector initialization \texttt{Y0\_init[:]}, we take
\begin{texteq} 
\begin{verbatim}
    Y0_init[:] = [0.1, 0.2, -0.33, -0.2, 0.5, -0.1]
\end{verbatim}
\end{texteq}
After 10 calls of the loss and gradient functions, the optimizer reaches the solution
\begin{texteq} \label{keplSolution}
\begin{verbatim}
    Y0[:] = [0.351, 0.706, -1.161, -0.238, 0.595, -0.12]
\end{verbatim}
\end{texteq}
which one checks corresponds to a bound state (\texttt{H\_kepl(y[:]=Y0[:])[] = -0.5 < 0}).
We can now evaluate~\texttt{NC\{\}},~\texttt{grad\_NC\{\}},~\texttt{hessian\_NC\{\}}
directly on the above solution. For the loss and the gradient, we obtain, respectively: 
\begin{texteq}
\begin{verbatim}
    NC(ydot{}=ydot_kepl, y_start[:]=Y0[:],
       t_final[]=T_final[])[]
    #  = 8.485708295471493e-19
    
    grad_NC(ydot{}=ydot_kepl, y_start[:]=Y0[:],
            t_final[]=T_final[])[:]
    #  = [0.0, 0.0, 0.0, 0.0, 0.0, 0.0]
\end{verbatim}
\end{texteq}
For the Hessian, we obtain the following eigenvalues:
\begin{texteq}
\begin{verbatim}
    &eigvals(hessian_NC(ydot{}=ydot_kepl, y_start[:]=Y0[:],
                        t_final[]=T_final[])[:,:])[:]
    #  =  [-2.39514e-07, -1.69088e-07, -8.688e-08,
    #      1.04679e-07, 5.13102e-07, 331.266786046988]
\end{verbatim}
\end{texteq}
In the above Hessian, we observe five flat directions.

These can be understood as follows: The dynamical symmetries of a bound
state of the Kepler problem form a $\mathfrak{so}(4)\simeq
\mathfrak{so}(3)\oplus\mathfrak{so}(3)$ Lie algebra, which is
six-dimensional. Picking the angular momentum algebra $(L_x, L_y, L_z)$ as
an obvious $\mathfrak{so}(3)$ subalgebra, we cannot find three conserved
quantities that commute with $\vec L$, so the angular momentum subalgebra
cannot be a regular subalgebra, it must be a diagonal. Then, the further
generators $A_i$ must transform nontrivially under the angular momentum
algebra. The only option here is that they transform as a vector, i.e. we
have $[\vec L,\vec L]\sim \vec L$, $[\vec L,\vec A]\sim \vec A$. This
leaves us with $[\vec A, \vec A]\sim\kappa\vec L$, where $\kappa$ provides
a sign that depends on the sign of $H$. For bound-states (``planetary
orbits''), we get $\mathfrak{so}(4)$ as a symmetry algebra, while for
scattering-states (``extrasolar comets''), we get $\mathfrak{so}(3,1)$,
with $\mathfrak{iso}(3)$ at-threshold (``parabolic-orbit'' bodies).  For
the Kepler problem, $\vec A$ is the Laplace-Runge-Lenz (`LRL') vector,
which for $k=m=1$ has the form $\vec p\times\vec L-\vec r/r$. Intuitively,
the symmetries associated with the components of the LRL vector do not turn
circular orbits into circular orbits, but rather change the eccentricity of
the orbit. Since they commute with the Hamiltonian, i.e. time-translation,
they leave energy, and also orbit-time, invariant. A useful geometrical
property here is codified by Kepler's third law -- Orbit-time is is a
function of the semimajor axis, i.e. `diameter' of the orbit. For a
circular orbit, there is a two-dimensional subspace of
eccentricity-inducing/adjusting deformations that leave the orbit-diameter
invariant, which can (for example) be parametrized by eccentricity and
angular position of the periapsis. The third component of the LRL vector --
the one parallel to angular momentum -- then acts trivially on the orbiting
body's state-of-motion. The situation parallels the observation that acting
with $SO(3)$ on a $3d$ unit-vector will only have two independent
$SO(3)$-generators change the vector, while the third (which generates
rotations whose axis is aligned with the given vector) acts trivially: In
the six-dimensional space of $SO(4)$-generators, we can find a
one-dimensional subspace of generators acting one acts trivially on any
given motion-state, and the other five give rise to flat directions in the
orbit-nonclosure Hessian.

In Figure~\ref{fig:orbitplots}(b), the full line shows the orbit whose initial configuration is the solution~\texttt{Y0[:]}.
We can use an SO(3)-rotation, as allowed by the symmetries of the system (conservation of angular momentum),
to rotate~\texttt{Y0[:]} into an initial configuration that evolves into an orbit in the (y,z) plane.
This is denoted in the figure by the dotted line.
The relevant rotation matrix is
\begin{texteq}
\begin{verbatim}
    R[:,:] = [[ 0.77615609,  0.40681284,  0.48175205],
              [-0.40681284,  0.90682312, -0.11034105],
              [-0.48175205, -0.11034105,  0.86933297]]
\end{verbatim}
\end{texteq}
which acts finitely on phase space in the usual way as
\begin{texteq} \label{keplSolutionDef}
\begin{verbatim}
    Y0[:] -> Y0_deformed[:] :=
      &concat(0, &es(R[:,:] @ i, j; Y0[:3] @ j -> i),
                 &es(R[:,:] @ i, j; Y0[3:] @ j -> i))
    #  = [1.11022302e-16, 6.25127871e-01, -1.25656979e+00,
    #     0.00000000e+00, 6.49574734e-01, -5.55264728e-02]
\end{verbatim}
\end{texteq}
The loss on this configuration is
\begin{texteq}
\begin{verbatim}
    NC(ydot{}=ydot_ho3, y_start[:]=Y0_deformed[:], t_final[]=T_final[])[]
    #  = 2.916662688870734e-16
\end{verbatim}
\end{texteq}

\subsection{The planar three-body-problem}

Our final example is the three-body-problem, restricted to motion in two
dimensions. We here use the common abbreviation `p3bp' for the planar
three-body problem. For simplicity, we consider units masses and unit
Newton's constant.  We first focus on the well-known figure-of-eight
orbit-closure solution~\cite{chenciner2000remarkable, moore1993braids}. We
subsequently perform a deformation of this solution along a flat direction
of its Hessian, which yields another orbit of the type discovered and
analyzed in~\cite{simo2002dynamical}.

A phase-space vector \texttt{y[:]} is now 12-dimensional: six positions
(two coordinates for each of the three masses) and six conjugate momenta
(again two coordinates per mass). Hamilton's equations for this system are:
\begin{texteq} \label{tbpsystem}
\begin{verbatim}
    ydot_p3bp = ^(
      y[:], 
      q[:,:] = &reshape(y[:6], 3, 2),
      dist[:,:,:] = q[:,+,:] - q[+,:,:], # dist[i,j,:] = q[i,:] - q[j,:]
      r_factor[:,:] = &pow(&es(dist[:,:,:] @ i, j, a;
                               dist[:,:,:] @ i, j, a -> i, j), -3/2))
       -> &concat(0,
            y[6:],  # dq/dt
            dist[1,0,:] * r_factor[1,0] + dist[2,0,:] * r_factor[2,0],
            dist[0,1,:] * r_factor[0,1] + dist[2,1,:] * r_factor[2,1],
            dist[0,2,:] * r_factor[0,2] + dist[1,2,:] * r_factor[1,2])[:]
\end{verbatim}
\end{texteq} 
where the intermediate quantity \texttt{r\_factor[i,j]},
\texttt{i,j = 0,1,2} computes the inverse cube of the
Euclidean separation between mass \texttt{i} and mass \texttt{j}.

As for the above two examples, our goal is to evaluate the gradient and
Hessian functions of the orbit-closure loss function on a closed orbit solution.
As discussed above, we aim for the figure-of-eight solution of~\cite{chenciner2000remarkable, moore1993braids}.
We quote the solution as given in the online three-body-problem gallery~\cite{BelgradeGallery, vsuvakov2013three},
recently compiled by the authors of~\cite{vsuvakov2013three} (their conventions align with our).
The initial vector is given as
\begin{texteq} \label{tbpInitVector}
\begin{verbatim}
    Y0_init[:] = [-1. ,  0. ,  1. ,
                  0. ,  0. ,  0. ,
                  0.347111,  0.532728,  0.347111,
                  0.532728, -0.694222, -1.065456]
\end{verbatim}
\end{texteq} 
The orbit time is quoted as
\begin{texteq}  \label{tbpTfinal}
\begin{verbatim}
T_final[] = 6.324449
\end{verbatim}
\end{texteq} 
and we fix it to this value for the remainder of this analysis.
Instead of evaluating the loss, gradient and Hessian directly on the configuration
in~Eq.~\eqref{tbpInitVector}, we first substitute it into a BFGS optimizer
(in the same way as we did in the Kepler problem). This allows us to flow
to a more accurate solution. After typically about 40 function calls,
the optimizer returns a solution such as
\begin{texteq} \label{tbpSolution}
\begin{verbatim}
    Y0[:] = [-9.99845589e-01, -5.69207692e-06,  9.99845620e-01,
             5.70200735e-06, -3.08148821e-08, -9.93042629e-09,
             3.47140692e-01,  5.32768073e-01, 3.47140612e-01,
             5.32768034e-01, -6.94281303e-01, -1.06553611e+00]
\end{verbatim}
\end{texteq}
with numerically-zero loss, zero gradient, and Hessian eigenvalues given as:
\begin{texteq} \label{tbpLossGradHess}
\begin{verbatim}
    NC(ydot{}=ydot_p3bp, y_start[:]=Y0[:],
       t_final[]=T_final[])[]
    #  = 4.459981866710229e-19
    
    grad_NC(ydot{}=ydot_p3bp, y_start[:]=Y0[:],
            t_final[]=T_final[])[:]
    #  = [6.98443059e-10,  8.23949237e-12, -7.23697523e-10,             
    #     8.16570974e-12, 2.52579189e-11, -1.64137619e-11,
    #     -1.97349744e-10, -2.79343635e-10, -1.97437565e-10,  
    #     -3.26573785e-10,  3.94809824e-10,  6.05863284e-10]
                                   
            
    &eigvals(hessian_NC(ydot{}=ydot_p3bp, y_start[:]=Y0[:],
                      t_final[]=T_final[])[:,:])[:]
    #  = [-2.9436638e-05, -2.214641e-06, 3.956369e-06,
    #     1.3811708e-05, 0.000595885249, 0.009097681599, 
    #     11.10411162849, 17.795125948157, 79.997311426776, 
    #     79.997322634127, 2626.009830021427, 10534.09893184725]
\end{verbatim}  
\end{texteq}
The Hessian eigenvalues in Eq.~\eqref{tbpLossGradHess} above suggest that the solution
Eq.~\eqref{tbpSolution} has four flat directions.
These may be identified as follows:
\begin{itemize}
\item time-translations, i.e.~shifting the initial state-vector along the orbit (1 flat direction). 
\item translations of the system in two spatial directions (2 flat directions).
\item SO(2)-rotation of the system in two spatial directions (1 flat direction).
\end{itemize}

We now discuss a deformation of this solution.
For the prior two systems (harmonic oscillator and Kepler problem) recall that the deformations given in Eq.~\eqref{hoscSolutionDef}, EQ.~\eqref{keplSolutionDef} of the original solutions Eq.~\eqref{keplSolution}, Eq.~\eqref{hoscSolution} were generated by the allowed symmetries of those solutions (a ~$GL(6)$ and~$SO(3)$ rotation respectively).
In this sense, the deformed solutions were equivalent, under symmetry, to the original solutions.

Here, we want to focus attention on a deformation that is \emph{not}
associated with a null-eigenvalue of the Hessian Clearly, such a
deformation cannot be symmetry of the problem, and the resulting
deformed state-of-motion vector will not correspond to a closed orbit.
That being said, we may pick an eigenvector with small associated
eigenvalue, so that we do less damage to orbit-closure than other
arbitrary deformations of similar magnitude would do.  Upon taking
this deformation, our approach is then to take the deformed vector as
initialization for a BFGS-type optimization step.

Now, from Eq.~\eqref{tbpLossGradHess}, we choose the eigenvector with
eigenvalue~$\approx 0.00059\ldots$.
The associated eigenvector is
\begin{texteq}
\begin{verbatim}
    deformation_evec[:] := &e_vec(hessian_NC(ydot{}=ydot_p3bp, y_start[:]=Y0[:],
                                             t_final[]=T_final[])[:,:])[:])[4]
                         = [-0.40539516 , -0.397770681 , 0.105262168 ,
                            0.289586383 , -0.312660227 , 0.0140926167 ,
                            0.241227186 , 0.406375493 , -0.0000000420869 ,
                           -0.0000000636894 , 0.0651632629 , 0.506912208]
\end{verbatim}  
\end{texteq}                                                
We apply a small shift along this direction to~\texttt{Y0[:]}:
\begin{texteq}
\begin{verbatim}
    Y0[:] -> Y0[:]_deformed_init[:] 
                 = Y0[:] - 0.02 * deformation_evec[:]                         
                 = [-0.9917376858 , 0.0079497215431 , 0.99774037664 ,
                    -0.0057860256527 , 0.0062531737251 , -0.0002818622644 ,
                    0.34231614828 , 0.52464056314 , 0.3471406128417 ,
                    0.5327680352738 , -0.695584568258 , -1.07567435416]          
\end{verbatim}  
\end{texteq}
where the magnitude~\texttt{0.02} of the deformation was chosen by hand\footnote{Note that, if the size of the shift along the non-symmetrical direction is taken to be too small, this minimization procedure will likely result in flowing back to the original solution.
However, if the shift-size is too large, the optimizer may not be able to flow to any solution at all.}.
We now pass~\texttt{Y0[:]\_deformed\_init} into the minimizer.
This leads to the new solution:
\begin{texteq} \label{tbpOrbitDeformation}
\begin{verbatim}
    Y0_deformed[:] = [-1.00149036,  0.00959776,  1.00543596, -0.01484088,
                      0.00831026, 0.00712495,  0.33774162,  0.53113049,
                      0.35292643,  0.53094484, -0.69066805, -1.06207533]
\end{verbatim}  
\end{texteq}
For this solution, the orbit-nonclosure loss is once again numerically-zero:
\begin{texteq} \label{tbpOrbitDeformationLoss}
\begin{verbatim}
    NC(ydot{}=ydot_p3bp, y_start[:]=Y0_deformed[:],
       t_final[]=T_final[])[]
                       = 3.834269651089502e-19
\end{verbatim}  
\end{texteq}
In Figure~\ref{fig:orbitplots}(c), we show the orbit corresponding to the original~\texttt{Y0[:]} solution Eq.~\eqref{tbpSolution} and the orbit corresponding to the deformed~\texttt{Y0\_deformed[:]} solution Eq.~\eqref{tbpOrbitDeformation}.
We note that, unlike for the original solution, the deformed solution gives trajectories where each mass moves on a~\textit{distinct} figure-of-eight curve. 
This class of solutions turns out to be known.
Indeed, they are referred to as the ``hyperbolic-elliptic companion orbits'' of the figure-of-eight orbit, having been discovered and analyzed in detail in work by Simo~\cite{simo2002dynamical}.

\newpage
\includepdf[pagecommand={\null\vfill\captionof{figure}{
\texttt{(a)} Orbit for the 3d harmonic oscillator dynamics in Eq.~\eqref{hoscsystem} with initial configuration (i.c.) in Eq.~\eqref{hoscSolution} and orbit time in Eq.~\eqref{hoscTfinal}. The dotted line is the deformed orbit with i.c.~in Eq.~\eqref{hoscSolutionDef}.
\texttt{(b)} Orbit for the 3d Kepler dynamics Eq.~\eqref{keplsystem} with i.c.~in Eq.~\eqref{keplSolution} and orbit time in Eq.~\eqref{keplTfinal}. The dotted line is the deformed orbit with i.c.~in Eq.~\eqref{keplSolutionDef}.
\texttt{(c)} Orbit for the 2d three-body dynamics Eq.~\eqref{tbpsystem} with i.c.~in Eq.~\eqref{tbpSolution} and orbit time Eq.~\eqref{tbpTfinal}. The dotted lines are the deformed orbit with i.c.~Eq.~\eqref{tbpOrbitDeformation}.
For each case, the orbit-non closure loss in Eq.~\eqref{Loss} for the solid-line orbit is included.
To visualize the velocity of a mass around its orbit, we include round markers which label its position at uniformly-spaced timesteps along the orbits.}\label{fig:orbitplots}},templatesize={145mm}{200mm},noautoscale=true,offset=-112 -50,scale=.85]{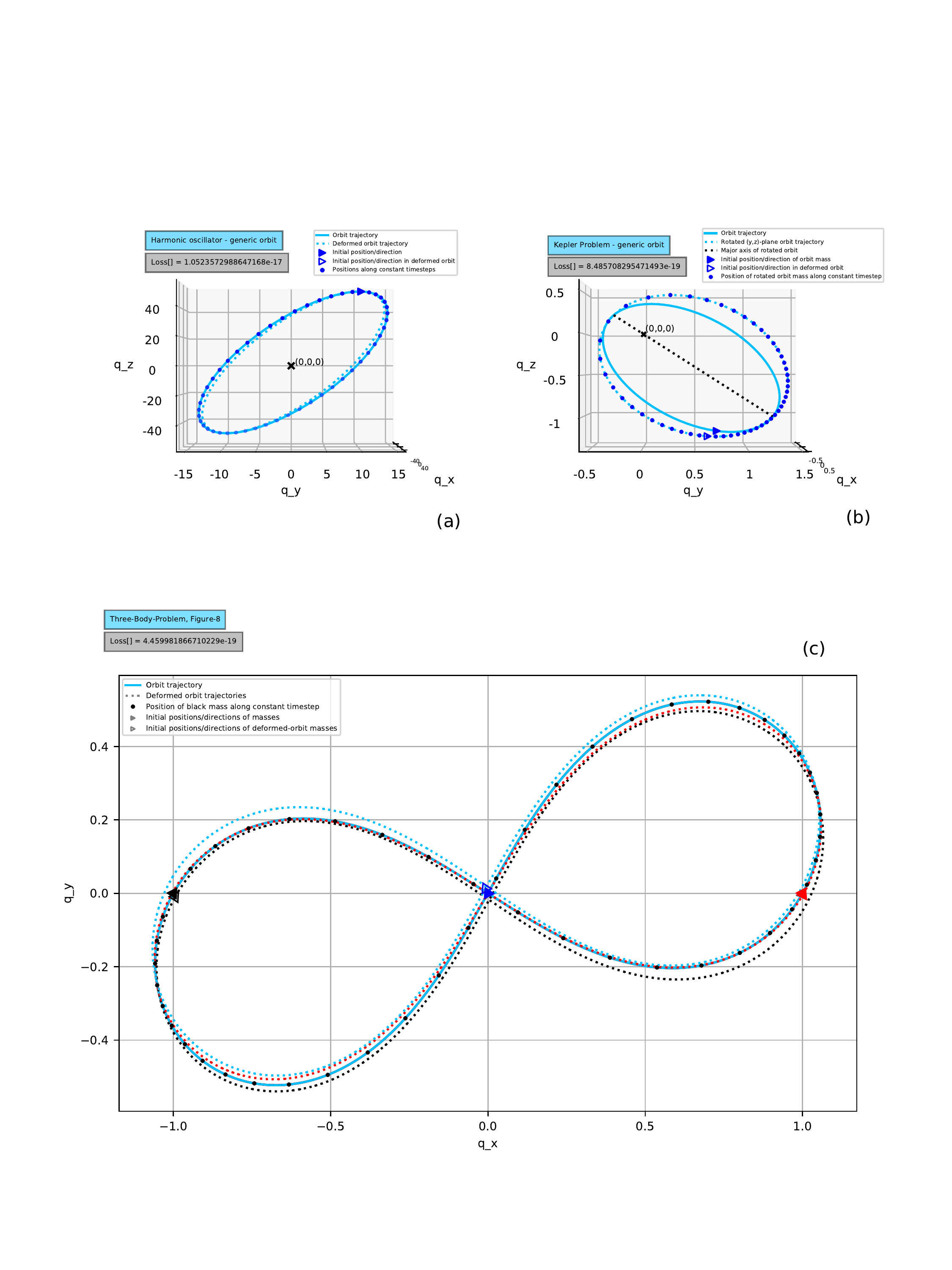}

\section{Conclusion and outlook}

The problem of backpropagating Hessians through ODEs involves quite a
few subtleties -- in general, the preferred choice of approach will
depend on characteristics of the problem, such as whether we want a
straightforward-to-understand scheme for validating more sophisticated
numerical code, whether or not it is feasible to employ some framework
that provides backpropagation capabilities, whether there are
disadvantages that are related to computing the Hessian all in one go
rather than row-by-row, and whether one is interested only in
backpropagating a Hessian around a minimum or not. Given that
sensitivity backpropagation allows us to compute the gradient of a
$\mathbb{R}^D\to\mathbb{R}$ scalar function that is defined in terms
of an algorithm with a constant-bound effort multiplier $b\le 8$, as
long as we can remember all intermediate quantities, and given also
that calculations that involve ODE integration admit forgetting
intermediate quantities if the ODE can be reversed in time, it is
unsurprising that an gradient of a scalar-valued function involving
ODE backpropagation with not remembered intermediate positions can be
computed numerically with effort multiplier $B=b+1$ (if using a fixed
time step ODE integrator whose computational effort characteristics
are easy to reason about), where the $+1$ is due to the need to redo
(and backpropagate) the velocity-field evaluation on the backward pass
in order to backpropagate also through that, for evaluating the
$\sigma_i F_{i,k}$-term. Then, using this approach for every component
of a gradient to obtain a Hessian, we can naturally expect to be able
to work this out with effort-multiplier no more than $B^2D$, not
remembering intermediate results for ODE-integration.

Each of the major approaches that have this time complexity in
principle has some drawbacks, such as: requiring some software
framework that cannot be used for the particular case at hand, not
automatically using a good starting point, requiring to solve an ODE
with $\le D^2$-dimensional state-vector, or, for the particular case
of hand-crafting backpropagation-of-backpropagation code, being
laborious and easily getting confusing and also not tolerating
constant-effort overhead for computing the rate-of-change function
(and its derivatives) well due to the large number of ODEs that need
to be solved. Pragmatically, one however observes that the
`differential programming' approach here is a very strong contender
for many real-world problems.

We addressed the point that `the wiring might get complicated' by
means of introducing a dedicated formalism. We hope that our thorough
discussion of various nontrivial aspects of the problem may provide
useful guidance to practitioners encountering such problems across a
broad range of fields.

The one aspect of this work that might appear most unfamiliar to a
numerical practitioner is the introduction of a formalism to help
manage the complexity that often arises when hand-backpropagating
code, including a fully-specified Chomsky type-2 grammar. Considering
the use case this formalism was designed for, it addresses the problem
that stringent equational reasoning should be possible for
`tensor-field like' quantities, as well as (potentially higher order)
functions involving such fields with complicated interdependencies,
\emph{and nothing beyond that}. This language in itself looks
promising for addressing a frequent complaint by newcomers to field
theory (including hydrodynamics) literature that `too many
dependencies are implicitly understood', and expressions hence often
are not really self-contained in the sense of a precise description of
the mathematical idea sans extra information of which dependencies
were suppressed. Many students first encounter this problem when
trying to make sense of the Euler-Lagrange approach to variational
calculus, in the form of ``$(d/dt)(\partial L/\partial\dot \vec
q)=\partial L/\partial \vec q$''.  One textbook that handles this
aspect with a similar emphasis as ours with respect to ``precisely
describing the wiring'' is~\cite{sussman2015structure}, and in a way,
the formalism presented there can be seen as an intellectual precursor
of (certainly an inspiration for) our formal language. Overall, there
are three major advantages from aligning the formalism with some
ingredients taken from the Python programming language, and NumPy:
insisting on named function-slots which may be both unbound or bound
not only allows good alignment with actual code implementations
(making the formalism suitable for code comments), but also does away
with the need of any `let'-syntax type syntactic sugar that many
functional programming languages identified as essentially necessary
to facilitate code clarity. Since our formalism has to support general
multilinear algebra, we deliberately do away with concepts such as
`matrices' and `vectors', which merely are $n=2$ and $n=1$ special
cases of $n$-index objects and use an uniform approach.

While we think (and hope) that our proposal to systematically describe
field-theoretic constructions of high complexity might be more
generally useful, we by no means consider the proposed language as
finalized. Overall, future amendments will however not go against the
grain of fundamental design principles. The most important principles
here are admitting and simplifying use in code-comments, to improve
alignment between expressions in research articles and in code,
and `staying minimalistic' in the sense of not introducing anything
that would introduce unnecessary special cases.

One aspect which we did not at all discuss in this work is in what
sense the ability to easily backpropagate Hessians might shed new
lights onto path integral based analysis, considering that in both
cases, we are often interested in exploring second-order spatial
dependencies of some objective function around a `classical'
trajectory.

\appendix

\section{Evaluation of $\sigma_i F_{i,k}$}
\label{app:velocity_jacobian}

Figure~\ref{fig:tf_jacobian} shows TensorFlow code that maps a
vector-valued TensorFlow function to a corresponding matrix-valued
function that evaluates the Jacobian, optionally admitting
autograph-transformation of the Python code for generating efficient
low-level code. Here, computation of the Jacobian will effectively use
backpropagation using, per intermediate quantity, one
sensitivity-accumulator per output-vector entry of the backpropagated
function. This means that for a $\mathbb{R}^D\to \mathbb{R}^N$ function~$F$,
the effort to compute the Jacobian is bounded by $bN$ times
the effort to evaluate~$F$ (with $b$ as in the main text).
In contradistinction, the code shown in figure~\ref{fig:tf_jv}
computes the gradient of a scalar-valued function $\sigma_j F_j$,
where only the $F_j$ are position-dependent.
This can be accomplished with effort-multiplier $b$.

\begin{figure}
  \label{fig:tf_jacobian}
{\small
\begin{lstlisting}[language=Python]
_DEFAULT_USE_AUTOGRAPH = True


def maybe_tf_function(use_autograph):
  """Returns tf.function if use_autograph is true, identity otherwise."""
  return tf.function if use_autograph else lambda f: f


def tf_jacobian(t_vec_func: Callable[[tf.Tensor], tf.Tensor],
                use_autograph=_DEFAULT_USE_AUTOGRAPH) -> Callable[[tf.Tensor],
                                                                  tf.Tensor]:
  """Maps a TF vector-valued function to its TF Jacobian-function."""
  # This is here only used to work out the hessian w.r.t. ODE initial-state
  # and final-state coordinates.
  # (!!!) Computing the costate-equation jacobian with this would be wasteful.
  @maybe_tf_function(use_autograph)
  def tf_j(t_xs):
    tape = tf.GradientTape()
    with tape:
      tape.watch(t_xs)
      v = t_vec_func(t_xs)
    ret = tape.jacobian(v, t_xs,
                        unconnected_gradients=tf.UnconnectedGradients.ZERO)
    return ret
  return tf_j
\end{lstlisting}
}
\begin{caption}
  {Mapping a vector-valued TensorFlow function to a function computing
   its Jacobian.  The
   \texttt{unconnected\_gradients=tf.UnconnectedGradients.ZERO}
   keyword argument ensures that TensorFlow produces a mathematically
   correct zero gradient, rather than a \texttt{None}-value if the
   graph for \texttt{v} is not connected to \texttt{t\_xs}.}
\label{fig:tf_jacobian}
\end{caption}
\end{figure}

\begin{figure}
{\small
  \begin{lstlisting}[language=Python]
def tf_jac_vec_v0(
    t_vec_func: Callable[[tf.Tensor], tf.Tensor],
    use_autograph=_DEFAULT_USE_AUTOGRAPH) -> Callable[[tf.Tensor, tf.Tensor],
                                                      tf.Tensor]:
  """Maps a TF vector-function F to a "(x, sx) -> sx_j Fj,k(x)" function."""
  @maybe_tf_function(use_autograph)
  def tf_j(t_xs, t_s_xs):
    tape = tf.GradientTape()
    with tape:
      tape.watch(t_xs)
      t_v = t_vec_func(t_xs)
    return tape.gradient(t_v, t_xs,
        output_gradients=[t_s_xs],
        unconnected_gradients=tf.UnconnectedGradients.ZERO)
  return tf_j


def tf_jac_vec_v1(
    t_vec_func: Callable[[tf.Tensor], tf.Tensor],
    use_autograph=_DEFAULT_USE_AUTOGRAPH) -> Callable[[tf.Tensor, tf.Tensor],
                                                      tf.Tensor]:
  """Maps a TF vector-function F to a "(x, sx) -> sx_j Fj,k(x)" function."""
  @maybe_tf_function(use_autograph)
  def tf_j(t_xs, t_s_xs):
    tape = tf.GradientTape()
    with tape:
      tape.watch(t_xs)
      t_v = tf.tensordot(t_s_xs, t_vec_func(t_xs), axes=1)
    return tape.gradient(t_v, t_xs,
        unconnected_gradients=tf.UnconnectedGradients.ZERO)
  return tf_j

\end{lstlisting}
}
\begin{caption}
  {Mapping a vector-valued TensorFlow function to the gradient of its
   scalar product with a given vector. This can be accomplished by
   either using the somewhat under-documented
   \codett{output\_gradients=} parameter of
   \codett{tf.GradientTape.gradient}, or alternatively by taking the
   partial gradient of a scalar product.}
\label{fig:tf_jv}
\end{caption}
\end{figure}

\section{A DSL for dependencies in tensor calculus}
\label{app:tcd_dsl}

This appendix describes the notational conventions used in this work
in detail. These conventions, which have been designed to simplify
equational reasoning about (multi-index) tensor calculus dependencies,
can be regarded as a domain specific language (DSL) -- henceforth
abbreviated `TCD-DSL'. They not only should be useful to adapt the
reasoning outlined in this work to minor variations such as explicitly
time-dependent velocity-functions, but may ultimately (perhaps after
some refinement of the grammar) be considered useful more widely for
generic vector/tensor calculus problems, not only at the interface of
symbolic analysis and numeric computation, but also whenever having a
precise specification of the wiring is important.

The incentive for this formal language comes from the need to simplify
keeping track of functional dependencies of tensorial quantities in
situations where the commonly used formalism that heavily makes use of
abuse-of-notation approaches becomes too unwieldy or ambiguous. The
grammar proposal presented here is to be seen as a starting point for
further discussion, and should not be considered as being cast in
stone.

While this article comes with a proper PEG grammar
specification~\cite{ford2004parsing} with which expressions can be
checked (using a parser generator) with respect to grammatical
well-formedness, no effort has yet been attempted to implement higher
level components of a compiler that can process computational
specifications following this grammar into various other forms,
including NumPy. Even having only the syntax analysis part of such a
compiler fleshed out clearly would be useful for discovering typos in
an automated fashion. Still, the structure of the grammar is such that
it should allow adding these components in a rather straightforward
manner. The TCD-DSL grammar in its present form is presented in
appendix~\ref{app:peggrammar}.

\subsection{Design desiderata}

The fundamental requirements that the language's design addresses are:

\begin{enumerate}
\item The language should be minimalistic, addressing precisely
  the unavoidable requirements, but with powerful primitives that
  eliminate the need for complicated structure.

\item The language should rely on tried-and-tested widely understood
  approaches for the sub-problems it needs to handle: some form of
  lambda calculus for functional dependencies and
  function-application, some form of tensor arithmetics to specify
  complicated tensor contractions and index reshufflings, and NumPy
  like syntax for common operations on multi-index arrays, such as
  reshaping, index-expansion, broadcasting, and slicing. This is in
  particular to simplify aligning formal expressions with code using a
  NumPy-like approach to handling multi-index data.

\item Expressions should be good to use in textual descriptions, such
  as research articles, but also admit being added verbatim as
  explanatory comments to code without requiring the reader to
  hallucinate the effects of typesetting instructions.
\end{enumerate}

The final requirement here asks for using a plain ASCII-based textual
representation, perhaps using indentation for emphasis.

We decided to use an amalgamation of lambda calculus and
python/lisp-style multi-parameter keyword arguments, with the
character `\codett{\^}' representing the `lambda` that indicates an
``anonymous function''.\footnote{In
ASCII-based typography, both the `backslash' and `hat' characters
historically have been used to represent the $\lambda$ of lambda
calculus. We here use the latter, since this blends in more nicely
both with text editors that tend to interpret a backslash as some form
of escaping for the subsequent opening parenthesis, and also
with \LaTeX{} typesetting.} The second requirement asks for
introducing some way to bind functions to names (treating functions as
`first-class citizens' in the value-universe) and to represent
function-evaluation. In a language like LISP or OCaml, it would not be
uncommon to have both a way to specify an anonymous function, such
as \codett{(lambda (x) (* x x))}, and on top of that some syntactic
sugar for introducing scoped variable bindings, such as \codett{(let
((x-squared (* x x))) \{body\})}, which is semantically equivalent to
an immediately-evaluated \codett{lambda}, such as \codett{((lambda
(x-squared) \{body\}) (* x x))}. If we require function parameters to
always be named, this allows us to kill three birds with one stone: We
get close alignment with a possible translation to code from
named-keyword parameters as available in Python (and also LISP), plus
a rather intuitive way to specify partial derivatives with respect to
parameters, and on top of that can meld \codett{let}
into \codett{lambda} with one simple convention.

\subsection{Concept overview}

The main concepts designed into this language are:

\begin{itemize}
\item
  Positionally-indexed multi-index numerical quantities.  Alongside
  function parameter names, the number of indices is generally
  indicated by square bracket indexing `\codett{[:, :, ...]}'  that
  would be valid Python syntax and a no-op when applied to a
  multi-index quantity of the intended shape, and not also valid for
  some quantity with a different number of indices: While a 3-index
  numpy array \codett{x} could be indexed as `\codett{x[:, :]}, this
  would also be valid indexing for a 2-index numpy array. This is not
  permitted here -- if a parameter is designated as \codett{x[:,:,:]},
  this clearly specifies that the quantity carries three indices.

\item
  Tensor contractions use an uniform approach that is based on
  generalized Einstein summation, via a special-function language
  primitive `\codett{\&es()}'. While this is slightly burdensome for
  simple scalar products and matrix/vector products, it follows the
  minimality principle that special cases of fundamental operations
  should not have special non-uniform syntax. As a consequence of
  this design, tensor index names (which are ubiquitous in field
  theory) can only occur to specify how to form and contract products,
  and outside such `\codett{\&es()}' expressions, indices on tensors
  are discriminated by index-order alone.

\item
  Immediately-bound function parameters. All function parameters are
  named, and unlike in pure lambda calculus (but like in most
  programming languages), functions can have multiple parameters.
  Borrowing some ideas from modern C++ constructor initialization
  syntax, named parameters can have a value assigned to them right at
  the time of their introduction, and this can refer to
  immediately-bound function parameters that were bound earlier
  (i.e. to the left of the current argument) on the same evaluation
  parameter-list. This means that we want to consider this expression,
  which describes a function with one open named scalar (zero-indices)
  parameter \codett{a} that gets mapped to \codett{a+a*a}
\begin{texteq}
\begin{verbatim}
    ^(a[], b[]=a[]*a[]) -> a[] + b[]
\end{verbatim}
\end{texteq}
as being equivalent to the expression:
\begin{texteq}
\begin{verbatim}
    ^(a[]) -> (^(a[], b[]) -> a[] + b[])(a[]=a[], b[]=a[]*a[])
\end{verbatim}
\end{texteq}

The corresponding rules for name-resolution are:

\begin{itemize}
\item If a name `\codett{x[...]}' is used in an expression, the name
  this is looked up by going through a sequence of fallbacks, where
  the first-found binding wins.

\item If the expression is the right hand
  side of a bound-parameter introduction, the name is looked up first
  among the earlier-introduced parameters on the same lambda
  parameter-list, where the latest (rightmost-in-text) such
  parameter-binding gets checked first.

\item Otherwise, the tightest lexically-surrounding lambda's parameter
  bindings, right-to-left, will be checked, then the lexically next
  outer lambda's bindings, etc. If no binding is found, the name is
  considered to be `unbound', so considered as coming from `global'
  context.
\end{itemize}

\item Partial derivatives. These are expressed by adding an index such
  as \codett{d Z[:]} to a positionally-indexed quantity. The semantic
  meaning of such an index is that, for every combination of the other
  indices present on the positionally-indexed quantity with this
  partial-derivative index removed, we use sensitivity backpropagation
  to determine the sensitivity of the corresponding coefficient on the
  coefficients of the input parameter named after the `d'. If a
  partial derivative is taken with respect to an immediately-bound
  multi-index quantity, the corresponding gradient represents the
  sensitivity (relative to $\varepsilon$, to first order) of the
  entries of the multi-index quantity without the partial-derivative
  slot with respect to $\varepsilon$-perturbations of the individual
  coordinates of the input. When taking partial derivatives with
  respect to a $k$-index input, this splices k~indices into the index
  list of the designated partial derivative (so, adding zero extra
  indices for a partial derivative with respect to a scalar).

  Hence,
\begin{texteq}
\begin{verbatim}
    ^(x[:], y[:]=x[:]) -> &es(x[:] @ a; y[:] @ a ->)[d y[:]]
\end{verbatim}
\end{texteq}
  is the derivative of the scalar product of a vector with itself
  with respect to varying the 2nd factor only. This is equivalent to
  the \codett{eps[i]}$\rightarrow 0$ limit of \codett{g} in:
\begin{texteq}  
\begin{verbatim}
    f = ^(x[:], y[:]) -> &es(x[:] @ a; y[:] @ a ->)[d y[:]]
    g = ^(x[:]) -> (f(x[:]=x[:], y[:]=x[:]+eps[:]) - f(x[:]=x[:], y[:]=x[:])) / eps[:]
\end{verbatim}
\end{texteq}
  using NumPy's broadcasting rules to divide by a vector.
  
\item Identifying a 0-parameter function with its value.
  We generally consider \codett{\^{}() -> T[:]} to be equivalent to \codett{T[:]} (and likewise for other numbers of indices on the object).

\item Functions can take functions as parameters. This is generally
  indicated by putting curly-braces indexing behind a slot-name.  The
  current form of the grammar makes no attempt to precisely represent
  the slot-name expectations for such function-parameters.

\item The grammar should also have a notion of ``whitespace'' that
  admits comments -- here in the form of ``hash character
  `\texttt{\#}' starts whitespace that ends at end-of-line''.
\end{itemize}

\subsection{Open questions}

It has been said that every formal language first and foremost
codifies what things its author did not properly understand. Major
pain points around the current language proposal which likely would
benefit from subsequent refinement include:

\begin{itemize}
  \item For function-parameters, no detail about their call signature
    is (or can be) made part of the parameter-specification. This
    is a trade off between brevity and precision.

  \item Ultimately, one likely would want to extend the grammar to
    admit tagging indices that correspond to a particular linear
    representation of some Lie group group and the representation,
    such as ``this is a SO(3, 1) co-vector index''.

  \item It is not yet clear what the best style decisions are with
    respect to proper use of technically unnecessary index-structure
    specifications. Is it visually clearer to mandate
  \codett{T(pos0[:]=v0[:], pos1[:]=v1[:])}, or should the terser form
  \codett{T(pos0[:]=v0, pos1[:]=v1)} be given preference?
    
\end{itemize}

\section{PEG grammar}
\label{app:peggrammar}

The PEG grammar for TCD-DSL is listed in this appendix. The grammar
rule comments provide additional clarification about specific rules.
This grammar can be directly processed by a parser generator such as
in particular the `lark' Python module~\footnote{Available at:
\url{https://github.com/lark-parser/lark}} to produce a working
parser. The code published alongside this article provides a
corresponding grammar-checker, which can also check expressions in
\LaTeX{} files.

\begin{scriptsize}
\begin{lstlisting}
// Valid input consists of zero or more statements.
?start: (signature_spec | statement)*

// A signature-specification.
signature_spec: (NAME ":")? LAMBDA lambda_args "=>" signature_rhs
?signature_rhs: (signature_spec | ("[" (COLON ("," COLON)*)? "]"))


// An individual statement is either an expression or a definition.
statement: expr | definition

// A definition is equivalent to a parameter binding.
// The only difference is that definitions are top-level, i.e. not scoped.
definition: parameter_binding

// An expression is a tensorial expression or a function-definition.
?expr: tensor_expr | fn_lambda

// Function, using "lambda notation": \(param1[:], param2[:]) -> {expr}
// Clarification: The right hand side expression can indeed
// be a function-definition, indicating a lexical closure.
fn_lambda: LAMBDA lambda_args "->" expr

// lambda arguments
lambda_args: "(" (fn_arg ("," fn_arg)*)? ")"

// Tensor expression
?tensor_expr: tensor_sum | tensor_optionally_indexed

?tensor_optionally_indexed: tensor_optionally_indexed_head bracketed_indexing*

?tensor_optionally_indexed_head: tensor_paren
  | tensor_special_fn
  | tensor_atom

// A parenthesized tensor expression may evaluate to a callable,
// so can optionally have call-arguments.
// This is valid:
//   "(\(x[:], y[:]) -> x[0]*y[0])(x[:]=y[:])"
// ("We take a function of two vectors, `x` and `y`, which computes the
//  product of the leading components, and bind the 2nd input to have
//  the same value as the first, leaving us with this function:
//  "\(x[:]) -> x[0] * x[0]".
?tensor_paren: /[+-]/? "(" expr ")" call_args*

// Function calls bind zero or more parameters.
// Interpretation: if only some slots are bound, the value of the
// call-expression is the function of the still-to-be-bound named
// parameters. A parameter-binding can refer to another earlier.
// parameter-binding on the same call.
call_args: "(" (parameter_binding ("," parameter_binding)*)? ")"

// A function argument is either a bound or an unbound slot.
fn_arg: slot_spec   // Unbound slot.
  | parameter_binding    // Bound slot.

// Q: Do we handle binding function-slots against functions properly?
parameter_binding: slot_spec "=" expr

// A tensor sum-expression can start with a sign, and is a sum of one
// or more tensor product-expressions.
?tensor_sum: /[+-]/? tensor_prod (/[+-]/ tensor_prod)*

// A tensor product-expression is a parenthesized or atomic
// tensor-expression, followed by zero or more such expressions
// separated by multiplication or division.
// The operator is made mandatory, since otherwise,
// telling apart a product expression such as "x[:](y[:]+z[:])"
// from a function-application such as "x[:](y[:]=z[:])"
// may be visually confusing.
?tensor_prod: tensor_optionally_indexed (/[*\/]/ tensor_optionally_indexed)*

// An atomic tensor expression is either a number or a name,
// optionally to be called, and optionally with indexing.
?tensor_atom: number | NAME call_args*

// Index-names may get extended later to also admit transporting
// information about the group and the representation.
index_name: NAME

index_name_seq: index_name ("," index_name)*

// ==============

// "[...]" indexing where individual indices may carry tags such as "@a".
bracketed_indexing: "[" (index ("," index)*)? "]"

// An index is either a Python "start:stop:range" type index, or "d slot".
// A '+', when used as an index, corresponds to `numpy.newindex`.
// NumPy's broadcasting rules are understood, except that broadcasting
// always requires both expressions to have the same number of indices,
// and if ranges for some particular index differ, then one of them
// must be a range-1 index.
index: DIFF WS diff_slot_spec
  | stride_part? (COLON stride_part?) ~ 0..2
  | "+"


// For ease of changing the definition of some tokens.
DIFF: "d"
LAMBDA: "^"
COLON: ":"

?stride_part: int_expr


// A slot specifier, such as `pos0[:]` in ODE(pos0[:]=y[:10])
// The only allowed forms are:
// - optional list of comma-separated COLONs in square brackets.
// - literal {} to indicate that we are passing in a function here.
//   TODO: When passing in a function, do we need to do something about slot-renaming?
slot_spec: NAME /\{\}/
  | diff_slot_spec


?diff_slot_spec: NAME ("[" (COLON ("," COLON)*)? "]")?


// Integer expressions - for indexing.
// Rules parallel the tensor_{paren|sum|prod} rules.
?int_num: /[+-]?[0-9]+/

?int_atom: int_num | /[+-]/? NAME

?int_expr: int_atom | int_paren | int_sum

?int_paren: /[+-]/? "(" (int_sum | int_paren) ")"

?int_sum: /[+-]/? int_prod (/[+-]/ int_prod)*

?int_prod: (int_paren | int_atom) (/[*\/]/? (int_paren | int_atom))*


// Tensor special-functions: concatenation, stacking, reshaping,
// and other operations.
// Also: Since all "field-functions" have all slots named,
// common unnamed-slots functions (like &exp) also will be of
// this type.
// These functions all start with a '&'.
tensor_special_fn: "&" any_special_fn

?any_special_fn: special_concat
 | special_stack
 | special_reshape
 | special_onehot
 | special_zeros
 | special_es
 | special_analytic

// Tensor-concatenation special function.
// Leading parameter needs to be a literal integer,
// the axis along which to concatenate.
special_concat: "concat" "(" int_num "," tensor_expr ("," tensor_expr)* ")"

// Stacking. Syntactically equivalent to concatenation.
special_stack: "stack" "(" int_num "," tensor_expr ("," tensor_expr)* ")"

// Reshaping. Indices are to be provided as a list of integers.
special_reshape: "reshape" "(" tensor_expr "," int_num ("," int_num)* ")"

// One-hot. Index and dimension can be int-expressions.
special_onehot: "onehot" "(" int_expr "," int_expr ")"

// Vector of zeros.
special_zeros: "zeros" "(" (NAME | int_num) ")"

// Generalized "Einstein Summation". This uses rather special syntax.
// - Product-factors are separated with semicolons,
// - Multi-index product-factors are expected to have all indices named
//   with "@ indexname1, indexname2, ...".
// - 0-index (scalar) tensors need not have index names.
// - result-indices are listed after '->'.
// Interpretation:
// - multi-index array indices are named according to the order in which
//   indices occur.
// - Every right hand side index combination addresses one coefficient
//   of the resulting tensor (in index order). The left hand side
//   indices that do not also occur on the right hand side
//   are being summed over. (This "generalized" Einstein summation is
//   group-theoretically unsound, but useful for formulas that are
//   ultimately not using group theory, but short-cuts.)
special_es: "es" "(" index_named ( ";" index_named ) * "->" index_names? ")"

// Single-argument special-functions,
special_analytic: special_analytic_1arg | special_analytic_2arg | special_matrix

// For mult-index tensorial quantities, these are understood to operate
// element-wise on the coordinates.
?special_analytic_1arg: ("exp" | "log") "(" tensor_expr ")"
?special_analytic_2arg: "pow" "(" tensor_expr "," tensor_expr ")"

// Eigenvalues/Eigenvectors, 2-index-to-1-index or 2-index-to-2-index:
?special_matrix: ("eigvals" | "eigvecs") "(" tensor_expr ")"


index_named: tensor_expr ("@" index_names?)?

index_names: index_name ("," index_name)*

NAME: /[A-Za-z][A-Za-z0-9_]*/

// TODO: This needs further refinement.
number: /[0-9]+/

// Whitespace -- will get ignored.
WS: /(\s|#.*)+/

%ignore WS
\end{lstlisting}
\end{scriptsize}

\bibliographystyle{unsrt}

\end{document}